\documentclass[journal,final]{IEEEtran}
\usepackage{cite} 
\usepackage{hyperref} 
\usepackage{xcolor}
\usepackage{multirow}
\usepackage{mathtools}
\usepackage{amsmath,amsfonts}
\usepackage{enumerate} 
\usepackage{algorithm}
\usepackage{algpseudocode}
\algrenewcommand\algorithmicrequire{\textbf{Input:}}
\algrenewcommand\algorithmicensure{\textbf{Output:}}
\newtheorem{theorem}{Theorem}
\newtheorem{definition}{Definition}
\newtheorem{lemma}{Lemma} 
\newtheorem{example}{Example} 
\usepackage{lscape}
\usepackage{adjustbox}
\usepackage{tabularx}
\usepackage{lipsum} 
\usepackage{array}

\begin{document}
\title{Helper and Equivalent Objectives: An Efficient Approach for Constrained Optimisation}
\author{Tao Xu and Jun He and Changjing Shang 
\thanks{Manuscript received xx xx xxxx}
\thanks{This work was partially supported by EPSRC under Grant No. EP/I009809/1. }
\thanks{Tao Xu and  Changjing Shang  are with the Department of Computer Science, Aberystwyth University, Aberystwyth SY23 3DB, U.K. E-mail:  \{tax2,cns\}@aber.ac.uk}
\thanks{Jun He is with the School of Science and Technology, Nottingham Trent University, Nottingham NG11 8NS, U.K. Email: jun.he@ntu.ac.uk (corresponding author)}
}
\maketitle

\begin{abstract}
Numerous multi-objective  evolutionary algorithms have been designed for constrained optimisation over past two decades. The idea behind these algorithms is to transform constrained optimisation problems into  multi-objective optimisation problems without any constraint, and then solve them. In this paper, we propose a new  multi-objective method for constrained optimisation, which works by converting a constrained optimisation problem into a problem with helper and equivalent objectives. An equivalent objective means that its optimal solution set is the same as that to the constrained problem but a helper objective does not. Then this multi-objective optimisation problem is decomposed into a group of sub-problems using the weighted sum approach. Weights are dynamically adjusted so that each  subproblem eventually tends to a problem with an equivalent objective.  We theoretically analyse the computation time of the helper and equivalent objective method on a hard problem called ``wide gap''. In a ``wide gap'' problem, an algorithm needs exponential time to cross between two fitness levels (a wide gap).  We prove that using helper and equivalent objectives can shorten the time of crossing the ``wide gap''.  We conduct a case study for validating our method. An algorithm with helper and equivalent objectives is implemented.  Experimental results show that its overall performance is ranked first when compared with other eight state-of-art evolutionary algorithms on   IEEE CEC2017 benchmarks in constrained optimisation.  
\end{abstract}

\begin{IEEEkeywords}
constrained optimisation,  constraint handling, evolutionary algorithms,  multi-objective optimisation, algorithm analysis, objective decomposition 
\end{IEEEkeywords}

\IEEEpeerreviewmaketitle

\section{Introduction}  
Optimisation problems in the real world usually are subject to some constraints. A single-objective constrained optimisation problem (COP) is formulated in a mathematical form as
\begin{align}
	\label{equCOP}
	\begin{array}{rll}
	\min  &f(\vec{x}), \quad \vec{x}=(x_1,\cdots,x_D)  \in \Omega,\\
	\mbox{subject to}&\left\{ 
	\begin{array}{lll}   
	g^I_i(\vec{x})\leq0,\quad    i=1,\cdots,q,\\
	g^E_i(\vec{x})=0,\quad i=1,\cdots, r,
	\end{array}\right.
	\end{array}
\end{align}
where  
$\Omega=\{ \vec{x}\mid  L_j \leq x_j\leq U_j, \; j =1, \cdots, D\} $  is a bounded domain in $\mathbb{R}^D$. $D$ is the dimension. $L_j$ and $U_j $ denote  lower and upper boundaries respectively.  
$g^I_i(\vec{x})\le 0$ is an inequality constraint and $g^E_i(\vec{x})=0$ is an equality constraint. A  feasible solution satisfies all constraints, and an infeasible solution violating at least one. The sets of optimal feasible solution(s), infeasible solutions and feasible solutions are denoted by $\Omega^*$, $\Omega_I$ and $\Omega_F$ respectively.  

Evolutionary algorithms (EAs) have been applied to solving COPs using different constraint handling methods, such as the penalty function, repairing infeasible solutions and multi-objective optimisation~\cite{michalewicz1996evolutionary,coello2002theoretical,mezura2011constraint,segura2016using}. 
A multi-objective method works by transforming a COP into a multi-objective optimisation problem  without inequality and equality constraints and then,   solving it by a  multi-objective EA. A popular implementation is to minimise the original objective function $f$ and  the degree of constraint violation $v$ simultaneously.
\begin{align}
\label{equBOP} 
  &\min \vec{f}(\vec{x}) = (f(\vec{x}),  v(\vec{x})), & \vec{x} \in \Omega. 
\end{align} 

The constraint violation degree in this paper is measured by the sum of each constraint violation degree.
\begin{eqnarray}\textstyle
 v(\vec{x})=\sum^q_{i=1} v^I_i(\vec{x})+\sum^{r}_{i=1} v^E_i (\vec{x}).
\end{eqnarray} 
$v^I_{i}(\vec{x})$ is the degree of violating the $i$th inequality constraint.
\begin{eqnarray}
&v^I_{i}(\vec{x})= \max \{ 0, g^I_i(\vec{x}) \}, & i=1, \cdots, q.
\end{eqnarray}
$v^E_i(\vec{x})$ is  the degree of violating the $i$th equal constraint.
\begin{eqnarray} 
&v^E_i(\vec{x})= \max \{0, \lvert g^E_i(\vec{x})\rvert-\epsilon\}, & i = 1, \cdots, r,
\end{eqnarray}
where $\epsilon$ is  a  user-defined tolerance allowed for the equality constraint. 

Multi-objective EAs  for constrained optimisation  have been proposed over past two decades. Many empirical studies have demonstrated the efficiency of the multi-objective method~\cite{segura2016using}.  Intuitively, the more objectives a problem has, the more complicated it is. Thus, this raises a question why the multi-objective method  could be superior to the single objective method. So far few theoretical analyses have been reported for answering this question.  

In fact, none of EAs in the latest IEEE CEC2017/18 constrained optimisation competitions adopted multi-objective optimisation~\cite{cec2017online}. The competition benchmark suite includes 50 and 100 dimensional functions. For a multi-objective optimisation problem, the higher dimension, the more complex Pareto optimal set. This raises another question whether multi-objective EAs are able to compete with the state-of-art single-objective EAs in the competition. 

The above questions motivate us to further study the multi-objective method for COPs.
Our work is inspired by  helper objectives~\cite{jensen2004helper}. 
The use of helper objectives has significantly improved the performance of EAs for solving some combinatorial optimisation  problems, such as job shop scheduling, travelling  salesman  and vertex covering~\cite{jensen2004helper,friedrich2009analyses}.  Our work is also inspired by objective decomposition, which was recently adopted in multi-objective EAs for COPs~\cite{xu2017new,zeng2017general,peng2018novel,wang2018decomposition}. Because the goal of COPs is to seek the optimal feasible solution(s) rather than a Pareto optimal set, decomposition-based multi-objective  EAs with  biased weights are flexible than those based on Pareto ranking. 

This paper presents a new equivalent and helper objectives method for COPs. A COP is converted into an optimisation problem consisting of equivalent and helper objectives but without any constraint. Here an equivalent objective means its optimal solution set is identical to $\Omega^*$, but a helper objective does not. Then this problem is solved by a decomposition-based multi-objective EA. 

Our research hypothesis is that the helper and equivalent objective method can outperform the single objective method on certain hard problems. We make both theoretical and empirical comparisons of these two methods. 

\begin{enumerate}
\item In theory, the ``wide gap'' problem~\cite{he2003towards,chen2010choosing} is regarded as a hard problem to EAs. We aim at proving   using helper and equivalent objectives can shorten the  hitting time of crossing such a ``wide gap''.  
    
\item A case study is conducted for validating our theory. We aim at designing an EA with helper and equivalent objectives and demonstrating that it can  outperform  EAs in  CEC2017/18  competitions.   
\end{enumerate}

This paper is a significant extension of  our two-page poster in GECCO2019~\cite{xu2019helper}. The  algorithm  described in the current paper is a slightly revised version of HECO-DE in~\cite{xu2019helper}. HECO-DE was ranked 1st in 2019 in IEEE CEC Competition on Constrained Real Parameter Optimization  when compared with other eight state-of-art EAs \cite{cec2017online}.

The paper is organised as follows:  Section~\ref{secFramework} is literature review. Section~\ref{secAnalysis} describes  the helper and equivalent objective method.   Section~\ref{secAnalysis} theoretically analyses this method. Section~\ref{secMOEA} conducts a case study. Section~\ref{secExperiments} reports experiments and results. Section~\ref{secConclusions} concludes the work.

\section{Literature Review}
\label{secReview}
Multi-objective EAs have been applied to COPs since 1990s~\cite{surry1997comoga,camponogara1997genetic}.  Segura et al.~\cite{segura2016using}  made a literature survey of the work up to 2016. Thus, this section focuses on reviewing most recent work. 
Following the taxonomy  in~\cite{mezura2008constrained,segura2016using}, a classification of these EAs is built upon  the type of objectives.
\begin{enumerate}
\item Scheme I  with two objectives,  the original objective $f$ and a degree of violating constraints $v$  \cite{zhou2003multi,cai2006multiobjective,wang2012combining,wang2012dynamic,xu2017new,peng2018novel,wang2018decomposition}.
\item Scheme II with many objectives,  he original objective $f$ and degrees of violating each constraint $v_i$ \cite{coello2002handling,li2017many}.
 
\item Scheme III with helper objective(s)     besides the original objective or the degree of constraint violation~\cite{deb2013bi,datta2016uniform,jiao2017dynamic,huang2019multiobjective}. For example, the penalty function forms helper objective. 
\end{enumerate}

The first scheme is the most widely used one so far. Ji et al. \cite{ji2017modified}  converted a  berth allocation problem with constraints into problem (\ref{equBOP}) and solved it by a modified non-dominated sorting genetic algorithm II. Ji et al. \cite{ji2018multiobjective} transformed a COP into problem (\ref{equBOP}) and solved it by a differential evolution (DE) algorithm. They combined multiobjective optimization with an $\epsilon$-constrained method. 

Recently, decomposition-based multi-objective EAs have applied to solving problem (\ref{equBOP}). Xu et al.~\cite{xu2017new} decomposed problem (\ref{equBOP}) into a tri-objective  problem  using the weighted sum method with static weights and solved it using a Pareto-ranking based DE algorithm.   Wang et al.~\cite{wang2018decomposition}  decomposed   problem~(\ref{equBOP}) using the weighted sum method into a number of subproblems with dynamical weights and solved these subproblems by DE.   Peng et al.~\cite{peng2018novel} decomposed  problem~(\ref{equBOP}) using the Chebyshev method. Weights are biased and adjusted dynamically for maintaining a balance between   convergence and  population diversity.

The second scheme converts a COP into a many-objective optimisation problem but is less used. Li et al.~\cite{li2017many} solved the many-objective optimization problem  by dynamical constraint handling.

The third scheme has an advantage of designing a helper objective. Zeng et al.~\cite{zeng2017general} designed a niche-count objective besides  the original objective and a constraint-violation objective and proposed an dynamic constrained multiobjective evolutionary algorithm (DCMOEA). The niche-count objective helps maintain  population diversity. They applied   three different  multiobjective  EAs (ranking-based,  decomposition-based,  and  hype-volume)  to the tri-objective optimisation problem. Jiao et al. \cite{jiao2017dynamic} converted a COP into a dynamical bi-objective optimisation problem consisting of the original objective and a niche-count objective.   Recently, these EAs with dynamic constrained multi-objectives were further improved by adding the feasible-ratio control technique  \cite{jiao2019feasible} and   a dynamic constraint boundary \cite{zeng2019constrained}.

The helper and equivalent objective method proposed in this paper belongs to the third scheme. One objective is designed as an equivalent objective. The equivalent objective has the same optimal set as that to the original COP. Helper objectives are also used to add more search directions. Under this framework, we have designed HECO-DE and HECO-PDE~\cite{huang2019multiobjective}.   HECO-PDE is an enhanced version of HECO-DE with principle component analysis. A multi-population implementation of HECO-DE is designed in~\cite{xu2019multi} which is suitable for parallel processing.

In order to speed up the convergence of EAs for COPs,   Deb and Datta \cite{deb2013bi} observed that the hybridisation of multi-objective EAs and local search can  reduce the number of fitness evaluations by one or more orders of magnitude. However, the current paper will not discuss the benefit of hybridisation but only focus on using helper and equivalent objectives.  

The  theoretical analysis of multi-objective EAs for constrained optimisation  is still rare and limited to combinatorial optimisation.  He et al.~\cite{he2014theoretical} proved that a multi-objective  EA with helper objectives is a 1/2-approximation algorithm for the knapsack problem. Recently, Neumann and  Sutton~\cite{neumann2018runtime} analysed the running time of a variant of Global Simple Evolutionary Multiobjective Optimizer  on the knapsack problem. Nevertheless, no general theoretical analysis  exists for the multi-objective EAs in continuous COPs.

\section{The Helper and Equivalent Objective Method}
\label{secFramework}
\subsection{Helper and Equivalent Objectives}
We start  from a problem existing in the classical bi-objective method for solving problem~(\ref{equBOP}). The Pareto optimal set to (\ref{equBOP}) is often significantly larger than $\Omega^*$. 
\begin{example}
\label{example1}
Consider the following COP. Its optimal solution is a single point $\Omega^*=\{0\}$. 
\begin{align*} 
\left\{
\begin{array}{rll}
\min &f(x)=x,  \qquad x \in [-1000,1000],\\
\mbox{subject to}  & g(x)=\sin(\frac{x \pi}{1200}) \ge 0.
\end{array}
\right.
\end{align*}
\end{example}

\begin{figure}[ht]
\begin{center}
  \includegraphics[height=3.3cm]{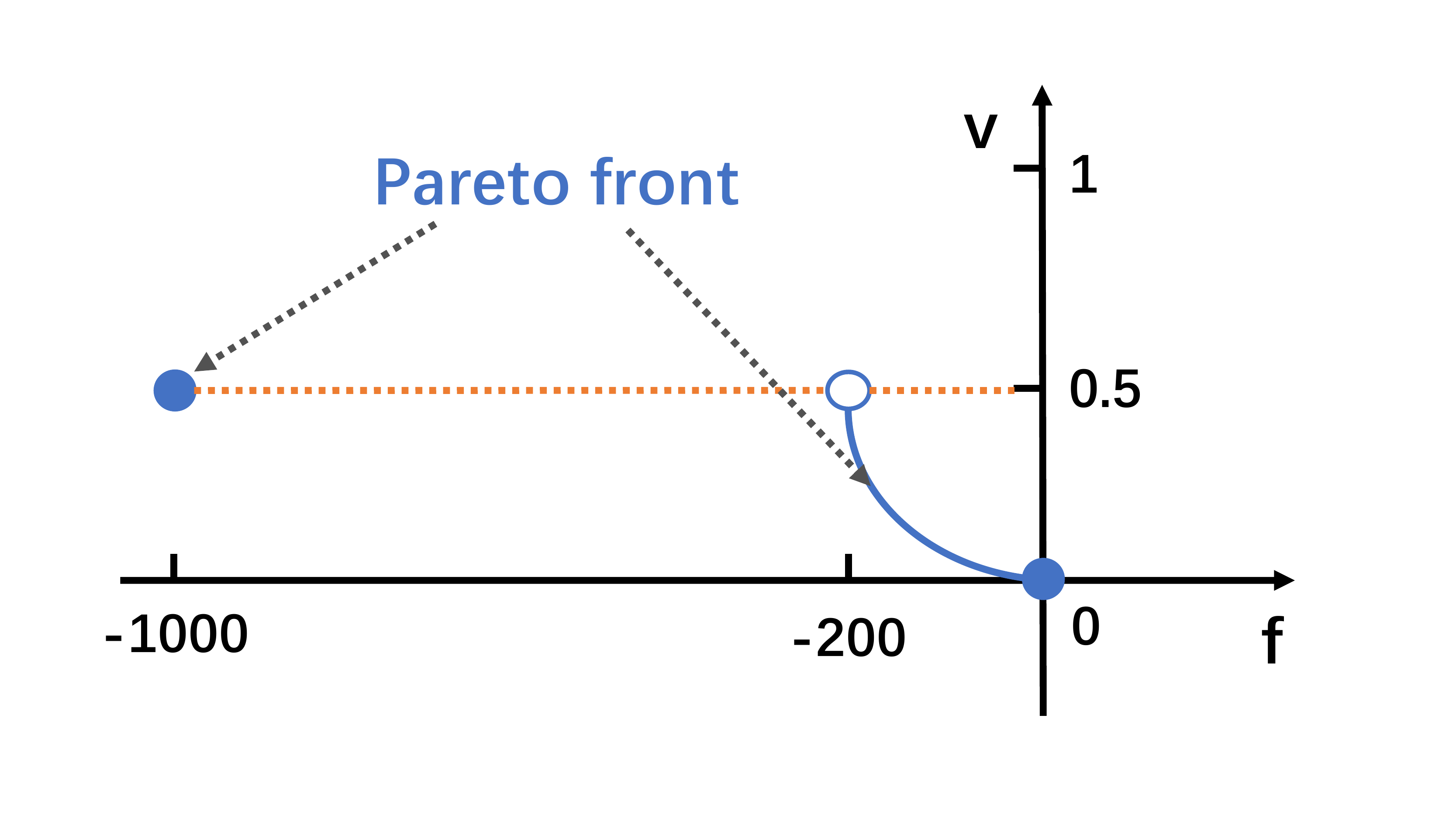}
\end{center} 
\caption{Pareto front.}
\label{fig1}
\end{figure}

The degree of constrain violation is
\begin{align} 
v(x)=\max \{0, -\sin({x \pi}/{1200})\}.
\end{align}
The Pareto optimal set to the bi-objective problem $\min (f,v)$ is $\{-1000\}\cup [-200,0]$, significantly larger than $\Omega^*$. The Pareto front is shown in Fig.~\ref{fig1}.

This example shows that using two objectives makes the problem more complicated. Thus, it is difficult to explain why  the multi-objective method is more efficient.

In order to develop a theory of understanding the multi-objective method for COPs, we introduce two concepts, equivalent and helper objectives. The term ``helper objective'' originates from~\cite{jensen2004helper}.  

\begin{definition} A scalar function $g(\vec{x})$ defined on $\Omega$ is called an  {equivalent objective function} with respect to the COP~(\ref{equCOP}) if it satisfies the condition:
\begin{equation}
    \arg \min \{f(\vec{x}); \vec{x} \in \Omega\} =\Omega^*.
\end{equation}
A scalar function $g(\vec{x})$ is called a   {helper objective function} if it does not satisfy the above condition.
\end{definition}

Equivalent functions can be obtained from  single objective methods for constrained optimisation. For example, a simple equivalent function is the death penalty function. 
Let $\Omega_F$ denote feasible solutions and $\Omega_I$ infeasible ones.  
\begin{align}
\label{equPenalty}
\min e(\vec{x})=
 \left\{
 \begin{array}{lll}
 f(\vec{x}), &\mbox{if } \vec{x} \in \Omega_F,\\
 +\infty , &\mbox{if } \vec{x} \in \Omega_I.
 \end{array}
 \right.
\end{align}

But the objective function $f$ is not an equivalent function unless all optimal solution(s) to $\min f$ are feasible. The constraint violation degree $v$ is not an equivalent function unless all feasible solutions are optimal. Hence, except particular COPs,  $\min (f, v)$  is a two helper objective problem.  

In practice, it is more convenient to construct an equivalent function  $e(\vec{x})$ which is defined on population $P$, rather than  $\Omega$. In this case, the definition of helper and equivalent functions is modified as follows. 

\begin{definition}
Given a population $P$ such that $\Omega^*\cap P \neq \emptyset$, a scalar function $g(\vec{x})$ defined on $P$ is called an  {equivalent objective function} with respect to the COP~(\ref{equCOP}) if it satisfies the following condition:
\begin{equation}
    \arg \min \{f(\vec{x}); \vec{x} \in \Omega \cap P\} =\Omega^* \cap P.
\end{equation}
A scalar function $g(\vec{x})$ defined on $P$ is called a   {helper objective function} if it does not satisfy the above condition. For a population $P$ such that $\Omega^*\cap P = \emptyset$,  we can not distinguish between equivalent and helper functions defined on the population. 
\end{definition}

An example is the superiority of feasibility rule~\cite{deb2000efficient} which is described as follows. Given a population $P$,
\begin{enumerate}
\item A feasible solution with a smaller $f$ value is better than one with a larger $f$ value;

\item  A feasible solution is better than an infeasible solution; 

\item  An infeasible solution  with smaller constraint violation is better than one with larger constraint violation.
\end{enumerate}  
 
The above rule leads to an equivalent function on $P$ as
\begin{align}
\label{equFeasibleRule}
e(\vec{x})=
 \left\{
 \begin{array}{lll}
 f(\vec{x}), &\mbox{if } \vec{x} \in \Omega_F\cap P,\\
 v (\vec{x})+f_F(P) , &\mbox{if } \vec{x} \in \Omega_I\cap P, 
 \end{array}
 \right.
\end{align}  
where   $
f_F(P)   =
 \max \{ f(\vec{x}), \vec{x} \in \Omega_F \cap P \}
$ if $\Omega_F \cap P \neq \emptyset$ or $f_F(P)=0$ otherwise.

\subsection{The Helper and Equivalent Objective Method}
Once an equivalent objective function is obtained, the COP~(\ref{equCOP})  can be converted to a single-objective optimisation problem  without any constraint. 
\begin{align}
\label{equSOP}
    \min  e(\vec{x}), \quad  \vec{x} \in  P.
\end{align}
In practice, an EA generates a population sequence $\{P_t; t=0,1, \cdots \}$ and $e(\vec{x})$ relies on population $P_t$.   

A single-objective EA  (SOCO) for problem (\ref{equSOP}) is described as follows.
 
\begin{algorithmic}[1]
\State  population $P_0 \leftarrow$ initialise a population of solutions; 
\For{$t=0, \cdots, T_{\max}$} 
\State    population $C_t \leftarrow$ generate a population of solutions from $P_{t}$ subject to a conditional probability $\Pr(C_{t}\mid P_t)$; 
\State    $P_{t+1}\leftarrow$ select  optimal solution(s) to $\min e(\vec{x}), \vec{x} \in P_t \cup C_t$; remove repeated solutions.    
\EndFor 
\end{algorithmic} 

$T_{\max}$ is the maximum number of generations. $\Pr(C_{t}\mid P_t)$ is a conditional probability determined by search operator(s). The population size $|P_t|$ is changeable so that $P_t$  is able to contain all found best solutions.

Besides the equivalent function $e(\vec{x})$, we add several helper functions $ h_i(\vec{x}), i=1, \cdots, k$, and then obtain a helper and equivalent objective optimisation problem on population $P$.  
\begin{align}
\label{equHEOP1}
\min \vec{f}(\vec{x}) =(e(\vec{x}),   h_{1}(\vec{x}), \cdots, h_k(\vec{x})), 
 &&
 \vec{x} \in  P.
\end{align} 
 
Furthermore, we decompose  problem~(\ref{equHEOP1}) into several single objective problem. Decomposition-based multi-objective EAs have been proven to be efficient in solving multiobjective optimisation problems~\cite{zhang2007moea,trivedi2017survey}.  
The decomposition method in the present work adopts the weighted sum approach, adding the helper objective onto the equivalent objective  such that
\begin{align}
\label{equWeightedSum}
\min \textstyle  w_0 e(\vec{x})  + \sum^k_{j=1} w_j h_j(\vec{x}), \quad  \vec{x} \in P,
\end{align} 
where  $w_j \ge 0$ are weights.  

Problem~(\ref{equHEOP1}) is transformed into $\lambda$ single-objective optimisation subproblems by assigning $\lambda$ tuples of weights $\vec{w}_i= (w_{0i}, w_{1i}, \cdots, w_{ki})$. 
\begin{align}
\label{equHEOP2} \textstyle
   \min f_i  =  w_{0i} e  + \sum^k_{j=1}w_{ji} h_j, \quad i=1, \cdots, \lambda.
\end{align} 
At least one $f_i$ is chosen to an equivalent objective function.  
We minimise all $f_i$ simultaneously.

Since the ranges of $e$ and $h$ might be significantly different, one of them may play a dominant role in the weighted sum. It is therefore, helpful to normalise the values of each function to $[0,1]$ so that none of them dominates others in the sum. The min-max normalisation method is adopted within a population $P$. Given a function $g(\vec{x})$, it is normalised to $[0,1]$.
\begin{align}
\label{equNorm}
    g(\vec{x}) \leftarrow \frac{g(\vec{x})-\max_{\vec{y} \in P} g(\vec{y} ) }{\max_{\vec{y} \in P} g(\vec{y})-\min_{\vec{y}} g(\vec{y})}.
\end{align} 

A helper and equivalent objective EA (HECO) for problem (\ref{equHEOP2}) is described as follows. 
 
\begin{algorithmic}[1]  
\State  population $P_0 \leftarrow$ initialise a population of solutions; 
\For{$t=0, \cdots, T_{\max}$}
\State  adjust weights; 
\State   population $C_t \leftarrow$ generate a population of solutions from $P_{t}$ subject to a conditional probability $\Pr(C_t \mid P_t)$; 
\State   $P_{t+1}\leftarrow$   select  optimal solution(s) to $\min f_i(\vec{x}), \vec{x} \in P_t \cup C_t$ for  $i=1, \cdots, \lambda$ where $f_i$ is calculated by formula (\ref{equHEOP2}); remove repeated solutions.  
\EndFor  
\end{algorithmic} 

HECO selects  optimal solution(s) to $\min  f_i(\vec{x}), \vec{x} \in P_t \cup C_t$ with respect to each function $f_i$ (called elitist selection), but it does not select all non-dominated solutions with respect to $(e, h_1, \cdots, h_k)$ (no  Pareto-based ranking).

Since our goal   is to find the optimal solution(s) to  
$\min e(\vec{x})$ but not to $\min h_i(\vec{x})$, it is not necessary to generate solutions evenly spreading on the Pareto front.   
Thus, the decomposition mechanism proposed herein differs from that employed in traditional decomposition-based multi-objective EAs~\cite{zhang2007moea}. The weights are chosen dynamically over generations $t$ so that each $f_i$ eventually converges to an equivalent objective function. Thus, the adjustment of weights  follows the   principle:  
\begin{align} 
\label{equWeights}
\lim_{t \to +\infty}   w_{0i,t} >0  \mbox{ and }\lim_{t \to +\infty} w_{ji,t}=0  \mbox{ for } j >0.
\end{align}

HECO has two  characteristics: 
\begin{enumerate}
\item   SOCO is one-dimension search along the direction $e$ in the objective space. HECO  is multi-dimensional search along  several directions $(e, h_1, \cdots, h_k)$.  $e$ is the main search direction for SOCO, while $h_1, \cdots, h_k$ are auxiliary directions added by HECO. Intuitively, if SOCO encounters a ``wide gap'' along the  direction $e$, HECO might bypass it through other auxiliary directions. This initiative discussion will be rigorously analysed later.

\item  The dynamically weighting ensures that at the beginning, HECO explores different directions $ e, h_1, \cdots, h_k$, while at the end, HECO exploits the direction $e$ for obtaining an optimal feasible solution.
\end{enumerate}
 
HECO is a general framework which covers many variant algorithm instances. Equivalent and helper functions can be constructed in a different way, such as (\ref{equPenalty}) and (\ref{equFeasibleRule}). Search operators can be chosen from evolutionary strategies, differential evolution, particle swarm optimisation and so on.

\subsection{Implicit Equivalent Objective}
Without the aid of an equivalent objective, a decomposition-based multi-objective EA  for COPs   faces a problem. The solution set found by the algorithm  is often larger than $\Omega^*$. This claim is shown through Example~1. We  assign $\lambda$ pairs of weights in objective decomposition: $(1,0), (w_i,1-w_i), (0,1)$ where $i=2, \cdots, \lambda-1$ and $w_i>0$ and obtain $\lambda$ subproblems with a bounded constraint $x \in [-1000,1000]$.
\begin{align*} 
\left\{
\begin{array}{ll}
\min f_1(x)=f,  \\
\min f_i(x)=  w_i f + (1-w_i) v ,  
\quad i=2, \cdots, \lambda-1,\\
\min f_{\lambda}(x)=  v.
\end{array}
\right.
\end{align*}

The  optimal solution to  $\min f$ is  $x=-1000$. The optimal solution to $\min f_i, i=2, \cdots, \lambda-1$ is infeasible. The optimal solution to $\min v$ is $[0,500]$.  The solution set to the $\lambda$ subproblems consists of infinite solutions, much larger than $\Omega^*=\{0\}$. Using dynamical adjustment of weights does not help here.

However, in practice, it is common to utilise the superiority of feasibility rule to select solutions. Using the rule, an infeasible solution such as $x=-1000$ is not selected. Among feasible solutions $x \in [0,500]$, only the minimal point $x=0$ is selected. But the superiority of feasibility rule is an equivalent objective (\ref{equFeasibleRule}), thus, many multi-objective EAs for COPs implicitly utilise an equivalent objective. Based on this argument, multi-objective EAs for COPs are classified into three types.
\begin{enumerate}
    \item Type I is to optimise helper objectives only;
    \item Type II is to optimise helper objectives but  select solutions by the superiority of feasibility rule (an implicit equivalent objective);
    \item  Type III is to explicitly optimise both helper and equivalent objectives.
\end{enumerate}

In this paper, the notation HECO refers to type III. It has some advantages:    an explicit equivalent objective is utilised and it can be designed more flexibly beyond the superiority of feasibility rule.  

\section{A Theoretical Analysis}
\label{secAnalysis} 

\subsection{Preliminary Definitions and  Lemma} 
Intuitively,  an equivalent objective ensures a primary search direction towards $\Omega^*$ and avoid an enlarged Pareto optimal set. Helper objectives provide auxiliary search directions. If there exists an obstacle like a ``wide gap'' on the primary direction, auxiliary directions can help bypass it. In theory, we aim at mathematically proving the conjecture: using helper and equivalent objectives can shorten the time of crossing the ``wide gap''.
First we introduce several preliminary definitions and a lemma.

For the sake of analysis, the search space $\Omega$ is regarded as a finite set. This simplification is  made due to two reasons. First, any computer can only represent a finite set of real numbers with a limited precision. Secondly,  population $P_t$ consists of  finite individuals (points). But the probability of $P_t$ at finite points always equals to $0$ in a continuous space. To handle this issue, we assume that possible values of $P_t$ are finite.     

Let $  \vec{f}(\vec{x})=(f_1(\vec{x}), \cdots, f_k(\vec{x}))$ be a scalar function ($k=1$) or a vector-valued function ($k> 1$). Consider a minimisation problem with bounded constraints:
\begin{align}
\label{equOP}
    \min \vec{f}(\vec{x}), && \vec{x} \in \Omega.
\end{align}
If $k=1$, it degenerates into a single-objective problem.
 
\begin{definition}  Given the optimisation problem (\ref{equOP}), 
$\vec{f}(\vec{x})$ is said to
\emph{dominate} $\vec{f}(\vec{y})$ (written as $\vec{f}(\vec{x}) \succ \vec{f}(\vec{y})$ ) if 
\begin{enumerate}
    \item $\forall i \in \{1, \cdots, k\}: f_i(\vec{x}) \le f_i(\vec{y})$;
    \item $\exists i \in \{1, \cdots, k\}:  f_i(\vec{x}) < f_i(\vec{y})$.
\end{enumerate}
If $k=1$, the two conditions  degenerate  into one  inequality $f(\vec{x})<f(\vec{y})$. 
\end{definition} 
 
Based on the domination relationship, the non-dominated set and Pareto optimal set are defined as follows.

\begin{definition}
A  set $S \subset S'$ is called  \emph{a non-dominated set in the set  $S'$} if and only if $\forall \vec{x} \in S$, $\forall \vec{y} \in S'$, $\vec{x}$ is not dominated by $\vec{y}$. A set $S$ is called \emph{a Pareto optimal set}  if and only if it is a non-dominated set in $\Omega$.  
\end{definition}

Given a target set, the hitting time is  the number of  generations for an EA to reach the set~\cite{he2017average}. The hitting time of an EA from one set to another is defined as follows. 

\begin{definition}  Let $\{P_t; t=0,1, \cdots\}$   be a population sequence of an EA.
Given two sets $S_1$ and $S_2$, the expected   hitting time of the EA from $S_1$ to  $S_2$ is defined by
\begin{align*}
    T(S_2\mid S_1):
    &\textstyle=\sum^{+\infty}_{t=0} \Pr(P_0 \subset  \overline{S_2},\cdots, P_t \subset  \overline{S_2}),
\end{align*} 
where the notation $\overline{S}$ denotes the complement set of $S$.  
\end{definition}

From the definition, it is straightforward to derive a lemma for comparing the hitting time  of two EAs.
\begin{lemma}
\label{lemma1} Let $\{P_t; t=0,1, \cdots\}$  and $\{P'_t; t=0,1, \cdots \}$ be two population sequences and $S_1$ and $S_2$ two sets such that $S_1 \cap S_2 =\emptyset$.  Let  $P_0=P'_0 = S_1$. If for any $t$,
\begin{align}
\label{equLemma1}
 \scriptstyle 
 \Pr(P_0 =S_1 \subset \overline{S_2},  \cdots, P_t \subset  \overline{S_2})  
 \ge \Pr(P'_0 =S_1 \subset \overline{S_2}, \cdots,   P'_t \subset \overline{S_2}), 
\end{align}
then $T(S_2\mid S_1)\ge T'(S_2 \mid S_1).$ Furthermore, if the inequality (\ref{equLemma1}) holds strictly for some $t$, then $T(S_2\mid S_1) > T'(S_2 \mid S_1).$
\end{lemma} 

This lemma provides a criterion to determine whether an EA has a shorter hitting time than another EA. The comparison is qualitative because no estimation of the hitting time is involved. For a quantitative comparison, it is necessary to utilise more advanced tools such as average drift analysis~\cite{he2017average}.  This will not be discussed in the current paper. 

\subsection{Fundamental Theorem}
Now we compare  SOCO  for  the single-objective problem~(\ref{equSOP}) 
and  HECO for the helper and equivalent objective problem~(\ref{equHEOP2}).
In order to make a fair comparison,  a natural premise is that both EAs use identical search operator(s).

The main purpose of using HECO is to tackle hard problems facing SOCO. 
Yet, what kind of problems are hard to SOCO? According to~\cite{he2003towards,chen2010choosing}, hard problems to EAs can be classified into two types: the ``wide gap'' problem and the ``long path'' problem. The concept of ``wide gap'' is established on fitness levels. In the helper and equivalent objective method, the equivalent function $e(\vec{x})$ plays the role of  ``fitness''. In constrained optimisation, function $f(\vec{x})$ is not suitable as ``fitness'' because the minimum value of $f$ might be obtained by an infeasible solution. 

The values of $e(\vec{x})$ are split into fitness levels: $FL_0< FL_1 <\cdots<FL_m$  and the search space $\Omega$ is split into disjoint level sets: $\Omega=\cup^m_{i=0} L_i$ where $L=\{\vec{x}; e(\vec{x})=FL\}$. Given a fitness level $FL$ and its corresponding point set $L$,
let $L^b$  denote points at better levels   $L^b:=\{ \vec{x} ;   e(\vec{x}) < FL\}$. 
A  ``wide gap'' between $L$ and $L^b$ is defined as follows.   

\begin{definition} Given an EA, we say  a  \emph{wide gap} existing between $L$ and $L^b$ if for a subset $A\subset L$,  the expected hitting time $T(L^b \mid  A \subset L)$ is an exponential function of the dimension $D$.  
\end{definition}

Several conditions are needed for mathematically comparing SOCO and HECO. 
Let $\{P_t;t=0,1, \cdots\}$ represent the population sequence from SOCO and $\{P'_t;t=0,1, \cdots\}$ from HECO.
Assume $P_0=P'_0$ are chosen from the fitness level $FL$. For  SOCO, thanks to elitist selection,  its offspring are either at the  level $FL$  or better fitness levels. 
For  HECO, because of selection on both equivalent and helper function directions, offspring  may include points from worse fitness levels  too. This  observation is summarised as a condition.

\textbf{Condition 1:}  Assume that $P_0 =P'_0 \subset L$. For SOCO, $P_t \subset L \cup L^b$ for ever. Provided that $P_t = X=(\vec{x}_1, \cdots, \vec{x}_m ) \subset L,$  
there is a one-to-many mapping from $P_t$ to $P'_t$ where $P'_t$ is    in the set $$
\scriptstyle Map(X)= \{X'=(\vec{x}_1, \cdots, \vec{x}_m, *)\mid * = \emptyset \mbox{ or }  *  \subset \overline{L \cup L^b} \}.$$  

The event of  $P_t= (\vec{x}_1, \cdots, \vec{x}_m) \subset L$  requires   $\vec{x}_1 \in L$, $\cdots,$ $ \vec{x}_m \in L$.   The probability of this event happening is larger than that of the event $P'_t= (\vec{x}_1, \cdots, \vec{x}_m, *) $ where $* = \emptyset \mbox{ or }  *  \subset \overline{L \cup L^b} \}$ because the latter event   requires    $\vec{x}_1 \in L$, $\cdots,$ $ \vec{x}_m \in L$ and also   $* \subset \overline{L \cup L^b}$.  This leads to the following conditions. 

\textbf{Condition 2:}   Let $P_{0}=P'_0=A \subset  L$. For any $t$, it holds 
\begin{align*} 
 &\scriptstyle\Pr(P_0 =A  \subset L, \cdots, P_t=Z \subset L) \\
 \ge &\scriptstyle\sum_{* \subset \overline{L^b} } \cdots \sum_{* \subset \overline{L^b} } \Pr(P'_0 =A' \subset \overline{L^b}, \cdots, P'_t=Z' \subset \overline{L^b}).
\end{align*} 

\textbf{Condition 3:}  For some $t$, the above inequality   is strict. 

Thanks to elitist selection and equivalent objective(s), Conditions~1 and 2 are always true. Condition 3 could be true, for example, if the transition probability from $*$ to $L^b$ is greater than 0.
Using the above conditions, we  prove a fundamental theorem of comparing HECO and SOCO. 
 
\begin{theorem}
\label{theorem1}
Consider  SOCO for  the single objective problem~(\ref{equSOP}) and HECO for  the helper and equivalent objective problem~(\ref{equHEOP2}) using elitist selection and identical search operator(s). 
Assume that SOCO faces a wide gap, that is, $T(L^b \mid A \subset L)$ is an exponential function of $D$ for a subset $A$. Let initial population $P_0=P'_0=A$.  Under Conditions~1 and~2, the expected hitting time  $T(L^b \mid A )\ge T'(L^b \mid A)$. Furthermore, under Condition~3, $T(L^b \mid A)>T'(L^b \mid A)$.
\end{theorem}

\begin{IEEEproof} 
From Conditions~1 and 2, it follows  for any $t$, 
\begin{align} 
&\scriptstyle\Pr(P_0 \subset  \overline{L^b},  \cdots, P_t\subset  \overline{L^b})  
=\scriptstyle  \sum_{A \subset L} \cdots \sum_{Z\subset L} \Pr(P_0=A ,  \cdots, P_t=Z ) \nonumber
\\\ge &\scriptstyle \Pr(P'_0\subset  \overline{L^b},  \cdots, P'_t\subset  \overline{L^b}) \nonumber\\
=&\scriptstyle \sum_{A \subset L} \cdots \sum_{Z\subset L}  \sum_{* \subset \overline{L^b} } \cdots \sum_{* \subset \overline{L^b} }  \Pr(P_0=A ,  \cdots, P_t=Z' ) . 
\end{align}
From Lemma~\ref{lemma1}, it is known  $T(L^b \mid A)\ge T'(L^b \mid A)$. The   second conclusion is drawn from Condition 3.
\end{IEEEproof}

Theorem~\ref{theorem1} proves  that the hitting time of HECO crossing a wide gap is not more than  SOCO under Conditions 1 and 2 (always true) and   shorter than SOCO under Condition 3 (sometimes true).   
In Conditions~2 and 3, the part $* \cdots *$ is a path of searching along helper directions and intuitively is regarded as a bypass over the wide gap. Theorem~\ref{theorem1} reveals if such a bypass exists, HECO may shorten the   hitting  time of  crossing the wide gap. Nevertheless, Theorem~\ref{theorem1} is inapplicable to the multi-helper objective method,  because the one-to-many mapping  in Condition~1  cannot be established.

\begin{example} Consider the  COP below,   
\begin{equation}
\label{equCaseCOP}
\left\{
\begin{array}{rll}
\min &f(x)=x,  \qquad x \in [-500,3000]\\
\mbox{subject to}  & g(x)=\sin(\frac{x \pi}{1000}) \ge 0.
\end{array}
\right.
\end{equation}
Its optimal solution is $x= 0$. The feasible region is  
$
    \Omega_F=[0,1000] \cup [2000, 3000].
$ The objective function $f(\vec{x})$ is not an equivalent function because its minimal point is $x=-500$, an infeasible solution.
\end{example}
 
First, we analyse a SOCO algorithm using elitist selection and the equivalent objective from the superiority of feasibility rule.   
\begin{align} 
\min e(x)=
 \left\{
 \begin{array}{lll}
 f(x), &\mbox{if } x \in \Omega_F,\\
 v (x)+3000 , &\mbox{if } x \in \Omega_I.
 \end{array}
 \right.
\end{align}  
where $v(x)=\max \{0, -\sin(\frac{x \pi}{1000})\}$. 

Mutation is     
$  y = x+U(-1,1),
$ 
where $x$ is the parent and $y$ its child. $U(-1,1)$ is a uniform  random number in $(-1,1)$.

Assume that SOCO starts at $L=\{2000\}$. Then  $L^b =[0,1000]$. Because of elitist selection, the EA cannot accept a worse solution. Then it cannot cross  the infeasible region $(1000,2000)$,  a wide gap to SOCO. Thus, $P_t \in L$ for ever. 

\begin{figure}[ht]
\begin{center}
  \includegraphics[height=3.3cm]{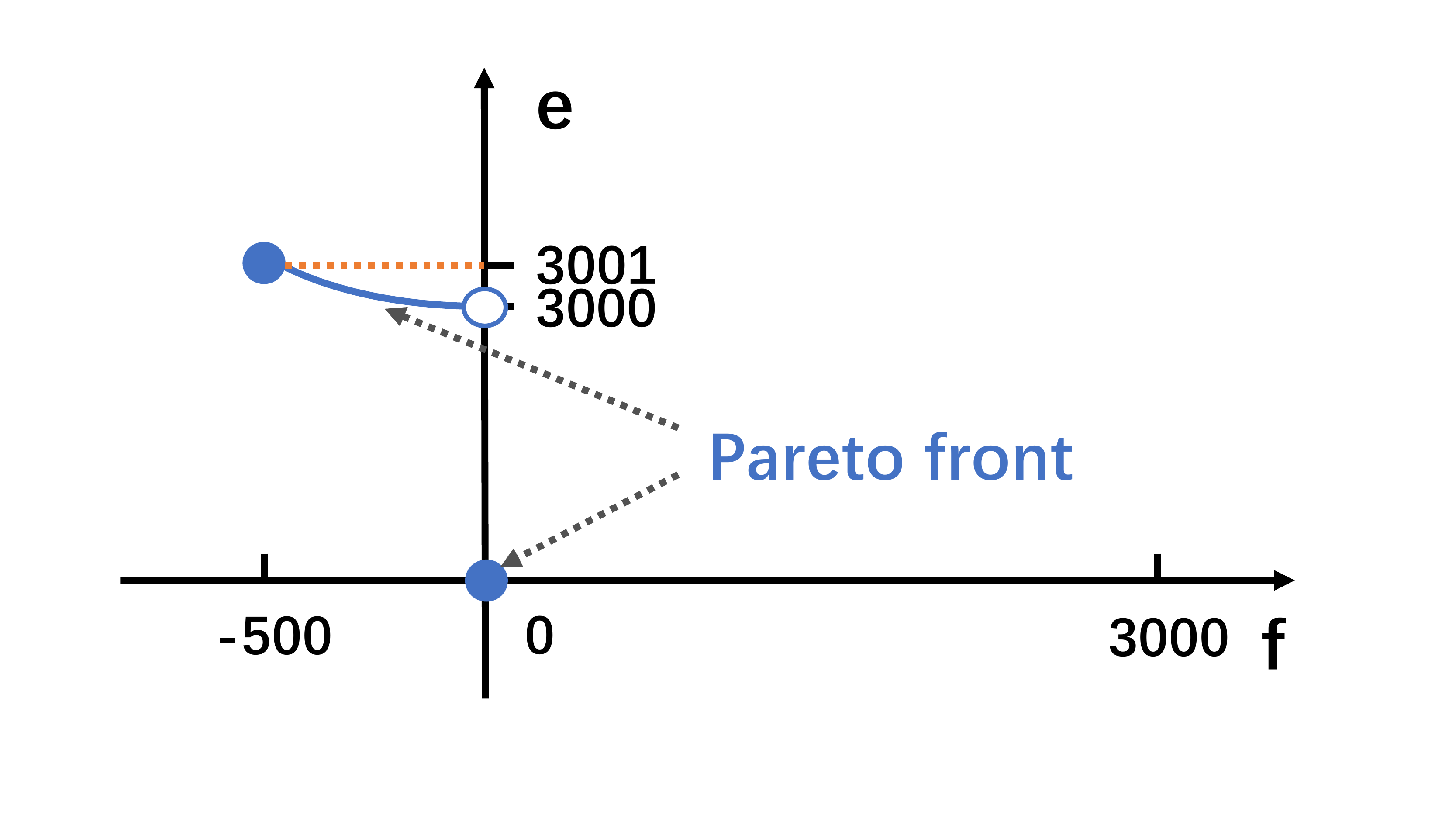}
\end{center} 
\caption{Pareto front  to the two-objective optimisation problem (22)}
\label{fig2}
\end{figure}

Secondly, we analyse a HECO algorithm employing elitist selection, identical mutation but two objectives.
\begin{align}
\label{equExample2}
\min \vec{f}(x)=(e(x),  f(x)), \qquad x \in [-500,3000].
\end{align} 
Its Pareto front is displayed in Fig.~\ref{fig2}.

\begin{figure}[ht]
\begin{center}
  \includegraphics[height=3.3cm]{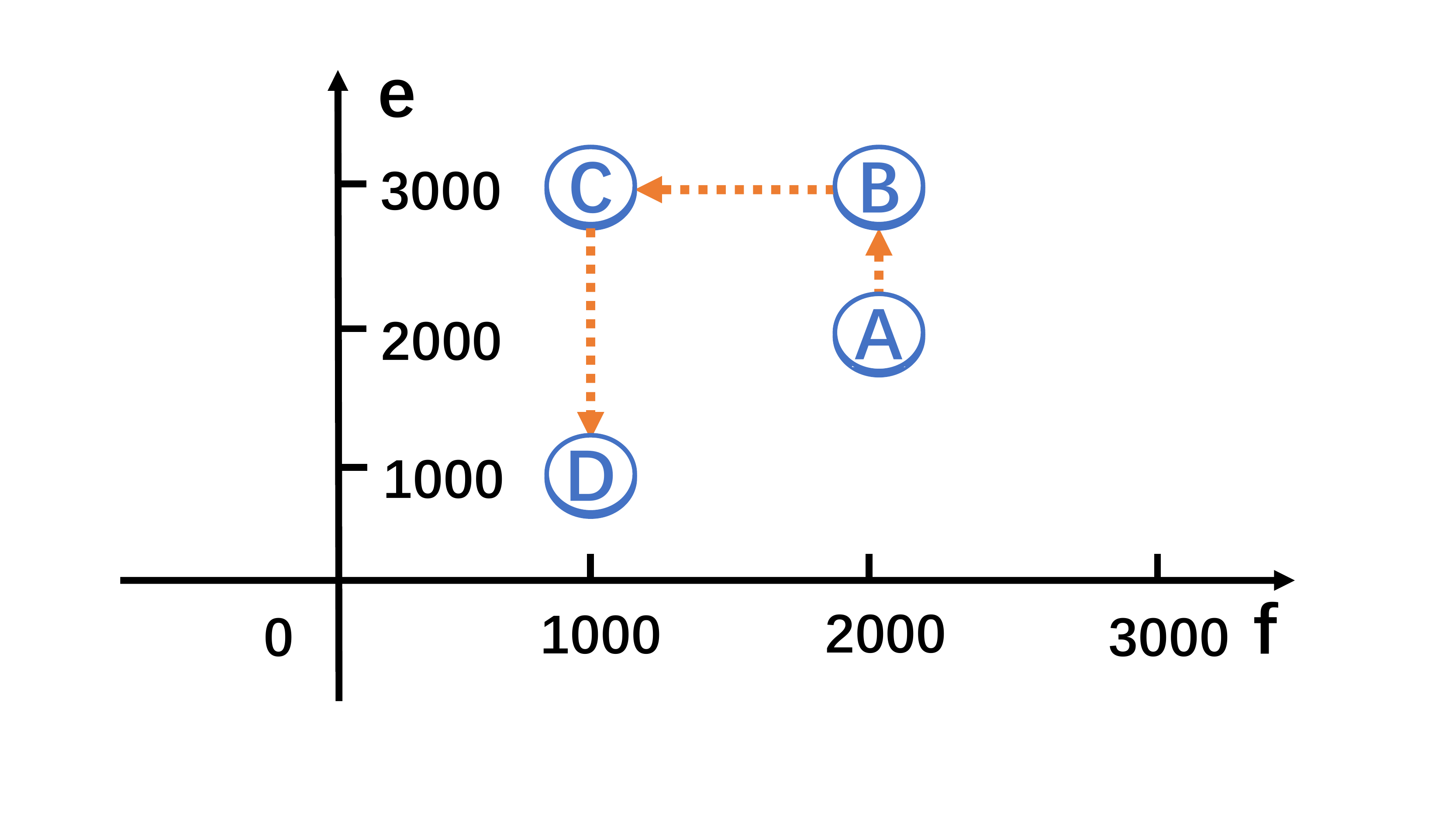}
\end{center} 
\caption{A bypass in objective space: $A(2000,2000) \to $ $ B(2000-\epsilon_1, 3000+\epsilon_2) \to$ $ C(1000-\epsilon_3, 3000+\epsilon_4) \to$ $ D(1000,1000)$ where $\epsilon_i \in (0,1)$ over the wide gap between fitness levels $e(x)=2000$ and $e(x)=1000$.}
\label{fig3}
\end{figure}

We assign two pairs of  weights: $\vec{w}_1=(1,0)$ and $\vec{w}_2=(0,1)$ on  $(e,f)$.
Assume that SOCO starts at $L=\{2000\}$.  For any $x \in P_t \cap [1000,2000]$, after mutation, some point $y$ such that $y <x-\frac{1}{2}$ is generated with a positive probability.  Since $f(y) < f(x)$, $y$ is selected to $P'_t$. Thus, $P'_t$ makes a downhill-search along the direction $f$. Repeating this procedure for 2000 generations, $P_t$ can reach  the set  $L^b=[0,1000]$ with a positive probability. This implies  for $t \ge 2000$,
 \begin{align*}  
\Pr(P'_0  \subset \overline{L^b}, \cdots, P'_t \subset  \overline{L^b})<1.
\end{align*}  
According to Theorem~\ref{theorem1},   $T'(L^b \mid L)< T(L^b \mid L)$. 
Fig.~\ref{fig3} visualises the bypass in the objective space.

\section{A Case Study}
\label{secMOEA} 
\subsection{Search Operators from LSHADE44} 
In order to validate our theory, we follow Occam's razor, that is to construct a HECO algorithm from a SOCO algorithm such that their search operators are identical but their objectives are different. No extra operation is added to HECO. For comparative purpose,  LSHADE44~\cite{polakova2017shade} is chosen as the SOCO algorithm  because it is ranked only 4th in the CEC2017/18 competition~\cite{cec2017online}. If  the constructed HECO algorithm  outperforms LSHADE44 and winer EAs in the  competition, then we have a good reason to claim the helper and equivalent objective method works. 

For the sake of a self-contained presentation,  search operators   in LSHADE44 are summarised as follows. 

LSHADE44 employs two mutation operators. The first one  is   current-to-pbest/1 mutation (see (6) in~\cite{zhang2009jade}).  Mutant point $\vec{u}_i$ is generated from  target point $\vec{x}_i$ by 
\begin{align}
    \label{equMutation1}
    \vec{u}_{i} = \vec{x}_{i} + F(\vec{x}_{pbest} - \vec{x}_{i}) + F(\vec{x}_{r_1} - \vec{x}_{r_2}),
\end{align} 
where  $\vec{x}_{pbest}$ is chosen at random from the top $100 p\%$ of population $P$ where $p\in (0,1)$.  $\vec{x}_{r_1}$  is  chosen at random from population $P$, while $\vec{x}_{r_2}$ at random from $P\cup A$ where $A$ represents an archive. Mutation factor $F \in (0,1)$.

The second mutation is randrl/1  mutation (see (3) in~\cite{kaelo2006numerical}).   
\begin{align}
    \label{equMutation2}
    \vec{u}_{i} = \vec{x}_{r_1} + F(\vec{x}_{r_2} - \vec{x}_{r_3}),\\
    \label{equMutation3}
    \vec{u}_{i} = \vec{x}_{r^*_1} + F(\vec{x}_{r^*_2} - \vec{x}_{r^*_3}).
\end{align}
In (\ref{equMutation2}), mutually distinct  $\vec{x}_{r_1}$, $\vec{x}_{r_2}$ and $\vec{x}_{r_3}$ are randomly chosen  from   population $P$. They are also different from $\vec{x}_{i}$.  
In (\ref{equMutation3}),  $\vec{x}_{r_1}$, $\vec{x}_{r_2}$ and $\vec{x}_{r_3}$ are chosen as that in (\ref{equMutation2}) but then 
are ranked. $\vec{x}_{r^*_1}$ denotes the best, while $\vec{x}_{r^*_2}$ and $\vec{x}_{r^*_3}$ denote the other two.  

LSHADE44 employs two crossover operators. The first one is binomial crossover (see (4) in~\cite{islam2011adaptive}). Trial point $\vec{y}_{i}$ is generated from target point $\vec{x}_{i}$ and mutant $\vec{u}_{i}$ by 
\begin{align}
\label{equCRbio}
    y_{i,j} = 
    \left\{
    \begin{array}{ll}
    u_{i,j},     &\mbox{if} \quad rand_j(0, 1) \leq CR \mbox{ or }  j = j_{rand}, \\
    x_{i,j},     &\mbox{otherwise},
    \end{array}
    \right.
\end{align}
where integer $j_{rand}$ is chosen at random from $[1, D]$. $rand_j(0, 1)$ is chosen at random from $(0,1)$. Crossover rate $CR\in[0, 1]$.
The second crossover is the exponential crossover (see (3) in \cite{zaharie2009influence}).

The combination of a mutation operator and a crossover operator forms a search strategy. Thus, four search strategies (combinations) can be produced. LSHADE44 employs a mechanism  of competition of strategies~\cite{tvrdik2006competitive,tvrdik2009adaptation} to create trial points. The $k$th strategy  is chosen subject to a probability $q_k$. All $q_k$ are initially set to the same value, i.e., $q_k = 1/4$. The $k${th} strategy is considered successful if a generated trial point $y$ is better than the original point $x$. The probability  $ q_k$ is adapted according to its success counts:
\begin{align}
\label{equStrat}
    q_k = \frac{n_k + n_0}{\sum_{i=1}^{4}(n_i + n_0)},
\end{align}
where $n_k$ is the   count of the $k${th} strategy’s successes, and $n_0>0$ is a constant.

LSHADE44 adapts parameters $F$ and $CR$ in each  strategy  based on previous successful values of $F$ and $CR$~\cite{polakova2017shade}. Each strategy has its own pair of memories $MF$ and $MC$ for saving $F$ and $CR$ values. The size of a historical memory is  $H$.

LSHADE44 uses an archive $A$ for the current-to-pbest/1 mutation~\cite{polakova2017shade}. The maximal size of archive $A$ is  set to $|A|_{\max}$.  At the beginning of search, the archive  is empty. During a generation, each point which is rewritten by its successful trial point is stored into the archive.  If the archive size exceeds the maximum size $|A|_{\max}$, then  $|A|-|A|_{\max}$ individuals are randomly removed from  $A$.

LSHADE44 takes a mechanism to linearly  decrease the population size~\cite{polakova2017shade,tanabe2014improving}. For population $P_t$, its size must equal to a required size $N_t$. Otherwise its size is reduced.  The  required initial size is set to $N_0$ and the finial size to $N_{T_{\max}}$. The required size at the $t$th generation is set by the formula:
\begin{align}
\label{equPopulationSize}
 \textstyle  N_{t} = round\left(N_{0} - \frac{t}{T_{\max}} (N_{0} - N_{T_{\max}})\right).
\end{align}
If $|P_{t}| >N_{t}$, then $|P_{t}|-N_{t}$ worst individuals
are deleted from the population.

\subsection{A New Equivalent Objective Function}
Two equivalent functions (\ref{equPenalty}) and (\ref{equFeasibleRule}) have been constructed from the death penalty method and  the superiority of feasibility rule respectively. However, measured by these functions, a feasible solution always dominates any infeasible one. To reduce the effect of such heavily imposed preference of feasible solutions,  we construct a new equivalent function. 

Let $\vec{x}^*_P$ be the best individual  in population $P$,  
\begin{align*}
    \vec{x}^*_P=
    \left\{
    \begin{array}{ll}
     \arg\min\{v(\vec{x}); \vec{x} \in P\},      & \mbox{if }   P\cap \Omega_F=\emptyset,  \\
    \arg\min\{f(\vec{x}); \vec{x} \in P\cap \Omega_F\},      & \mbox{if }   P \cap \Omega_F\neq \emptyset.
    \end{array}
    \right.
\end{align*}
For each $\vec{x} \in P$, $\tilde{e}(\vec{x})$ denotes the fitness difference between $f(\vec{x})$ and $f(x^*_P)$.
\begin{align}
\label{equ:neweq}
\tilde{e}(\vec{x}) =|f(\vec{x})-f(\vec{x}^*_P)|
\end{align}

$\tilde{e}$ itself is not an equivalent function because in some problems, the fitness  of an infeasible solution  is equal to $f(\vec{x}^*_P)$ too. An equivalent function  on population $P$ is defined as
\begin{equation}
 \label{equEquivalentFunction}
    e(\vec{x})= 
    w_1 \tilde{e}(\vec{x}) +w_2 v(\vec{x}), 
\end{equation} 
where  $w_1, w_2>0$ are weights, which are used to control the contribution of  $\tilde{e}$ and   $v$ to the equivalent function $e$.
The number of such equivalent functions is infinite because  $w_1 \in (0,+\infty), w_2 \in (0, +\infty)$.  

\begin{theorem}
Function $e(\vec{x})$ given by (\ref{equEquivalentFunction}) is an equivalent objective function for any weights $w_1>0, w_2>0$. 
\end{theorem}

\begin{IEEEproof}
Given any $P$ satisfying $\Omega^*\cap P \neq \emptyset$,  we have $\min \{e(\vec{x}); \vec{x} \in P\}=0$.  On one hand, for any $\vec{x} \in \Omega^* \cap P$, $\hat{e}(\vec{x})=0$ and $v(\vec{x})=0$, then $e(\vec{x})=0$. On the other hand, for $\vec{x} \in P$ such that $e(\vec{x})=0$,  it holds $v(\vec{x})=0$, then $\vec{x} \in \Omega^*$.
\end{IEEEproof}

\begin{figure}[ht]
     \centering
     \includegraphics[width=70mm]{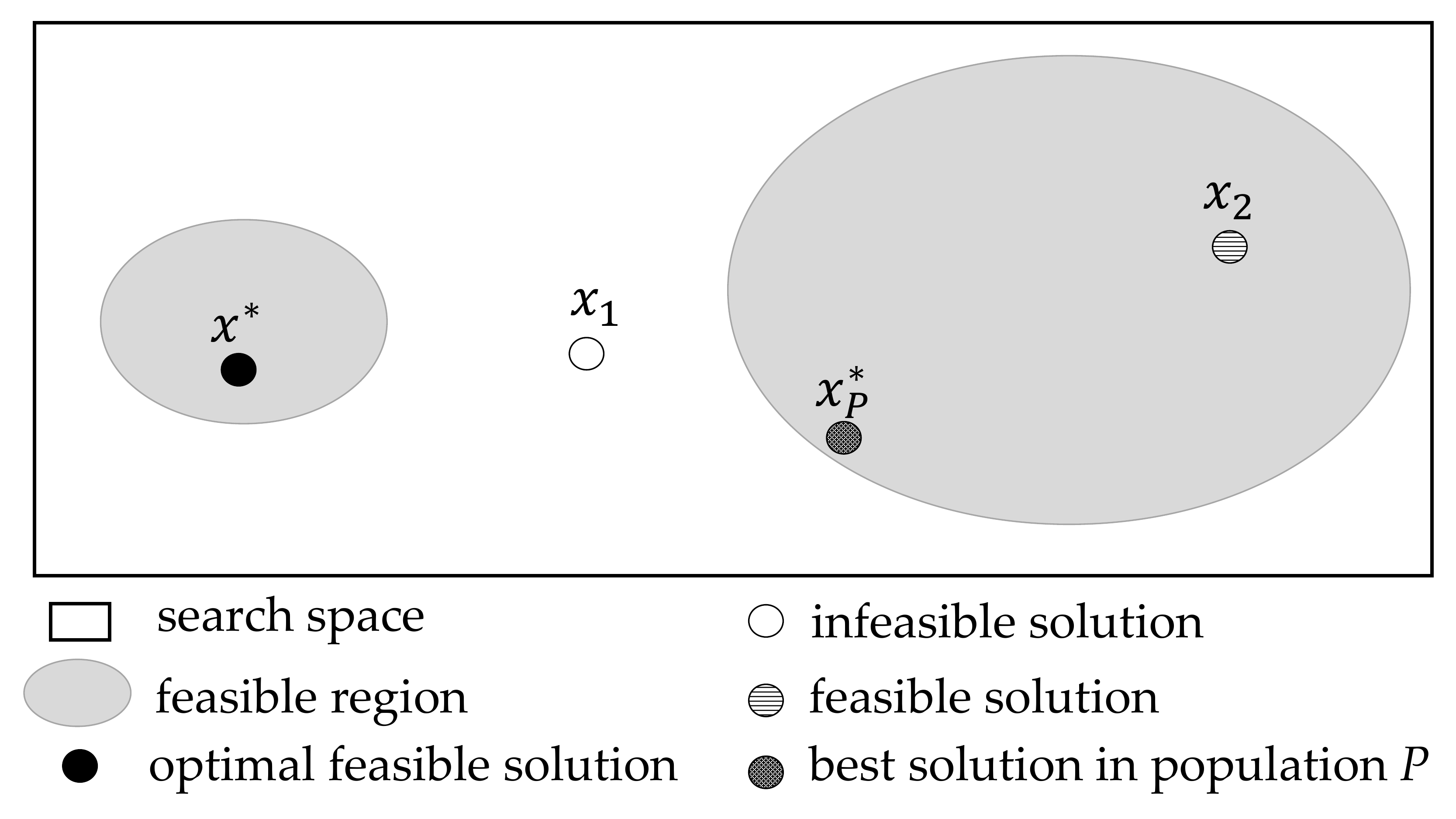}
     \caption{There exist two feasible regions. An infeasible $\vec{x}_1$ satisfying (\ref{equCompara}) is better than $\vec{x}_2$ under the equivalent objective function $e$. This may help population $P=(\vec{x}_1, \vec{x}_2, \vec{x}^*_P)$ move from the right feasible region to the left feasible region in which the optimal feasible solution $\vec{x}^*$ locates.}
     \label{fig4}
 \end{figure}

If two solutions $\vec{x}_1$ (infeasible) and  $\vec{x}_2$ (feasible) in population $P$ satisfy  
\begin{align}
\label{equCompara} 
w_1 |f(\vec{x}_1)-f(\vec{x}^*_P) |+w_2 v(\vec{x}_1) <  w_1 |f(\vec{x}_2)-f(\vec{x}^*_P)|,
\end{align}
then  under the equivalent objective function $e$, infeasible $\vec{x}_1$ is better than feasible $\vec{x}_2$.  This feature may help search the infeasible region. For example, in Fig.~\ref{fig4}, assume that $f(\vec{x}_1)-f(\vec{x}^*_P)=0$ and $f(\vec{x}_2) -f(\vec{x}^*_P) =1$, $v(\vec{x})_1=0.5$  and $w_1=w_2$. Then we have $e(\vec{x}_1) =0.5 e(\vec{x}_2)$.  
Starting from $\vec{x}_1$, it is much easier to reach the left feasible region in which the optimal feasible solution $\vec{x}^*_P$ locates.

We choose $f$ as a helper function and then obtain a problem with helper and equivalent objectives.  
\begin{align}
\label{equHEOP3}
\min \vec{f}(\vec{x}) =(e(\vec{x}),   f(\vec{x})), 
 &&
 \vec{x} \in P, 
\end{align} The problem is decomposed into $\lambda$  single objective subproblems through the weighted sum method: for $ 
   i=1, \cdots, \lambda,$ 
\begin{align}
\label{equHEOP4}
\begin{array}{rr}
     &   \min f_i(\vec{x})=  w_{1i} \tilde{e}(\vec{x}) + w_{2i} v(\vec{x})+w_{3i} f(\vec{x}). 
\end{array}
\end{align} 
An extra term $\tilde{e}$ is added besides the original objective function $f$ and constraint violation degree $v$.  


\subsection{ A New multi-objective EA for Constrained Optimisation}
A HECO algorithm is designed which reuses search operators from LSHADE44~\cite{polakova2017shade}. We call it HECO-DE because it is built upon HECO and DE. Different from the single-objective method LSHADE44, HECO-DE has three new multi-objective features: helper and equivalent objectives, objective decomposition and dynamical adjustment of weights.
The procedure of HECO-DE is described in detail as below. 

\begin{algorithmic}[1]
    \State Initialise algorithm parameters, including  the required initial population sizes $N_{0}$ and final size  $N_{T_{\max}}$, the maximum number of fitness evaluations $FES_{\max}$, circle memories for parameters $F$ and $CR$, the size of historical memories  $H$; initial probabilities  $q_k$ of four strategies, and external archive $A$;  
    \State Set the counter of fitness evaluations $FES$  to $0$, and the counter of generations $t$  to $0$;
    \State Randomly generate  $N_{0}$ solutions and form an initial population $P_0$;
    \State  Evaluate the value of $f(\vec{x})$  and  $v(\vec{x})$ for each  $\vec{x} \in P_0$;
    \State Increase counter $FES$ by $N_{0}$;  
    \While {$FES \le FES_{\max}$ (or $t \le T_{\max}$)}
        \State Adjust weights in objective decomposition.  
        \State Assign sets $S_{F}$ and $S_{{CR}}$ to $\emptyset$ for each strategy. The sets are used to preserve successful values of $F$ and $CR$  for each search strategy respectively.  The set $C$ (used for saving children population)  is also set to $\emptyset$. 
        \State Randomly select $\lambda$ individuals (denoted by $Q$) from $P$ and then denote the rest individuals  $P \setminus Q$ by $P'$;  

        \For{$x_i$ in $Q$, $i = 1,\dots,\lambda$}
            \State  Select one strategy (say $k$) with  probability $q_k$ and generate mutation factor $F$ and crossover rate $CR$ from respective circle memories;  
            \State  
            Generate a trail point $\vec{y_i}$ by applying the selected strategy;  
            \State {Evaluate the value of $f(\vec{y_i})$ and $v(\vec{y_i})$;}
            \State Add $\vec{y_i}$ to  subpopulation $Q$, resulting in an enlarged subpopulation $Q'$; 
            \State Normalise $\tilde{e}(\vec{x})$, $f(\vec{x})$ and $v(\vec{x})$ for each individual $\vec{x}$ in   $Q'$.
            \State  Calculate $f_i$ value for $\vec{x}_i$  and $\vec{y}_i$ according to formula~(\ref{equHEOP4}).
            \If {$f_i(\vec{y_i}) < f_i(\vec{x_i})$}
               \State {Add $\vec{y}_i$ into children $C$ and $\vec{x}_i$ into archive $A$;} 
                \State {Save values of $F$ and $CR$ into respective sets $S_{F}$ and $S_{CR}$ and increase respective success count;}
            \EndIf
        \EndFor
        \State  Update circle memories $M_{F}$ and $M_{CR}$  using respective sets $S_F$ and $S_{CR}$ for each strategy (see its detail in LSHADE44~\cite{polakova2017shade});
        \State Merge subpopulation $P'$ (not involved in mutation and crossover) and children $C$ and form new population $P$;
        \State {Calculate the required population size $N_{t}$;}
        \If {$N_{t} < |P|$}
            \State{Randomly delete $|P|-N_{t} $ individuals from $P$;}
        \EndIf
         \State Calculate the required archive size $|A|_{\max} =4 N_{t}$; 
        \If {$|A| > |A|_{\max}$ }
            \State{Randomly delete $|A| - |A|_{\max}$ individuals from archive $A$};
        \EndIf
        \State Increase  counter   $FES$  by $\lambda$ and  counter   $t$  by 1;
    \EndWhile  
\end{algorithmic}  

There are several major differences between HECO-DE and LSHADE44 which are listed as below. 

\textit{Lines 12}: in HECO-DE, mutation is applied to subpopulation $Q$, rather than the whole population $P$. Thus,    current-to-pbest/1 mutation and randr1/1 mutation must be modified because the ranking of individuals is   restricted to subpopulation $Q$.  Given target $x_i$ and subpopulation $Q$,   $x_{Qbest}$ is chosen to be the individual in   $Q$ with the lowest value of $f_i(\vec{x})$.  Hence, current-to-pbest/1 mutation (\ref{equMutation1}) is modified as 
\begin{align}
    \label{equMutation1'}
    \vec{u}_{i} = \vec{x}_{i} + F_k(\vec{x}_{Qbest} - \vec{x}_{i}) + F_k(\vec{x}_{r_1} - \vec{x}_{r_2}),
\end{align}  
This new mutation is called  current-to-Qbest/1 mutation. For randr1/1 mutation~(\ref{equMutation3}), $\vec{x}_{r_1}$, $\vec{x}_{r_2}$ and $\vec{x}_{r_3}$ are not compared but just randomly selected from subpopulation $Q$. Thus it returns to the original rand/1  mutation~(\ref{equMutation2}).

\textit{Lines 12 and 16}: ranking  individuals  is used in both mutation (\ref{equMutation1}) and calculation of the equivalent function (\ref{equEquivalentFunction}). Because ranking is restricted within subpopulation $Q$ and its size $\lambda$ is  a small constant, the time complexity of ranking is a constant.  This is different from LSHADE44 in which individuals in the whole population $P$ are ranked. Its time complexity is a function of dimension $D$.

\textit{Lines 17-20}: if $f_i(\vec{y_i}) < f_i(\vec{x_i})$,   then $\vec{y}_i$ is accepted and added into children population $C$. HECO-DE minimises $\lambda$  functions $f_i$  simultaneously.  In \textit{Line~7}, the weights on each $f_i$ are dynamically adjusted (detail in Subsection~\ref{sec:Weights}).
This is the most important difference from LSHADE44.

Since $\lambda$ is a small constant, the number of operations in HECO-DE  is only changed by a constant when compared with LSHADE44. Thus, the time complexity of HECO-DE  in each generation is the same as  LSHADE44~\cite{polakova2017shade}.

\subsection{A New Mechanism of Dynamical Adjustment of Weights}  
\label{sec:Weights}
We propose a special mechanism for dynamically adjusting weights. Function $f_i$ in   subproblem~(\ref{equHEOP4}) is a weighted sum of helper and equivalent functions:
\begin{align} 
\label{eq:fi}
 f_i(\vec{x})=  w_{1i} \tilde{e}(\vec{x})  +w_{2i} v(\vec{x}) + w_{3i} f(\vec{x}), 
\end{align} 
where  $w_{1i}, w_{2i}, w_{3i}$ are the weights on functions $\tilde{e}, v$ and $f$ respectively.   
Weights  are adjusted according to the following principle: each   $f_i$ converges to an equivalent function. Thus,
$$\lim_{t \to +\infty}  w_{1i,t} >0,   \lim_{t \to +\infty}  w_{2i,t}  >0,  
  \lim_{t \to +\infty}  w_{3i,t} =0.$$   

In HECO-DE, weights are designed to linearly increase (for $w_1, w_2$) or decrease (for $w_3$) over $t$ and also linearly increase (for $w_1, w_2$) or decrease (for $w_3$) over $i$.  In more detail, weights are  given by  
\begin{align}
    w_{1i,t} &= \frac{t}{T_{\max}}\cdot\frac{i}{\lambda}, \label{equW1}\\
    w_{2i,t} &=  \frac{t}{T_{\max}}\cdot\frac{i}{\lambda} + \gamma, \label{equW2}
\\
    w_{3i,t} &= \left(1 - \frac{t}{T_{\max}}\right)\cdot\left(1 - \frac{i}{\lambda}\right), \label{equW3}
\end{align} 
where $\lambda$ is the number of subproblems.   $T_{\max}$  is the maximal number of generations. $\gamma \in (0,1)$ is a bias constant which is linked to the number of constraints. The more constraints, the larger $\gamma$ and $w_2$.  

\begin{figure}[ht]
\centering
  \includegraphics[width=0.45\textwidth]{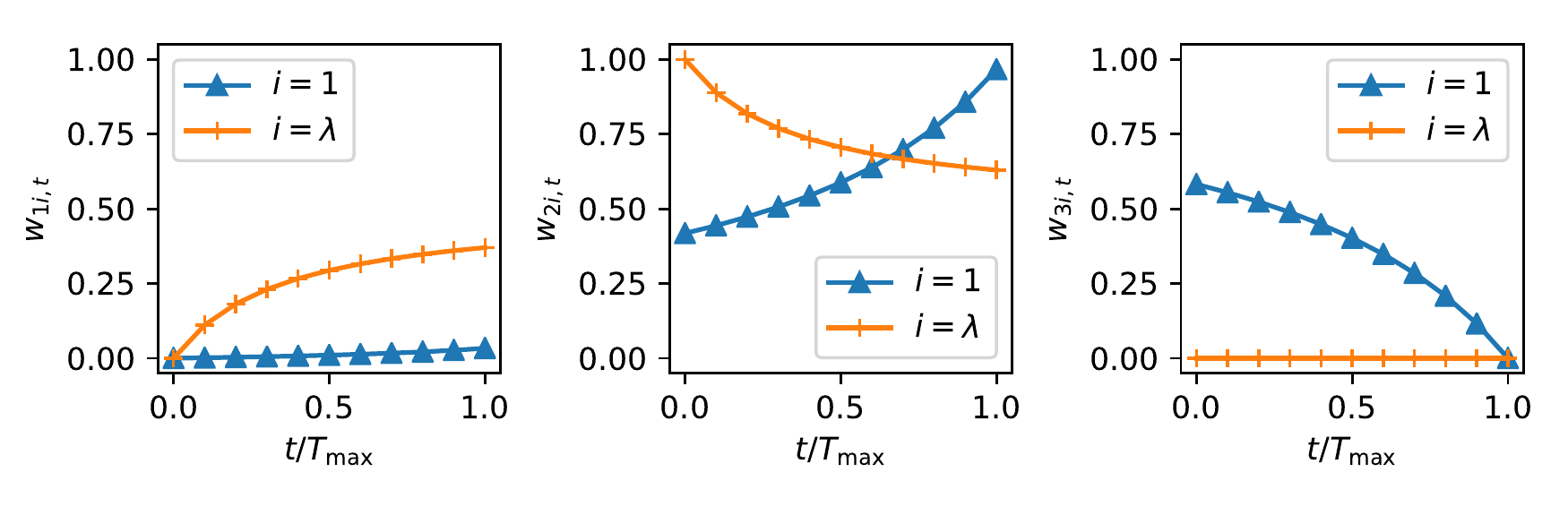}
 \caption{The change of weights for $1$st and $\lambda$th individuals on  CEC2006 benchmark functions.  $\gamma =0.7$.}
  \label{fig5} 
\end{figure}

\begin{figure}[ht]
\centering
  \includegraphics[width=0.45\textwidth]{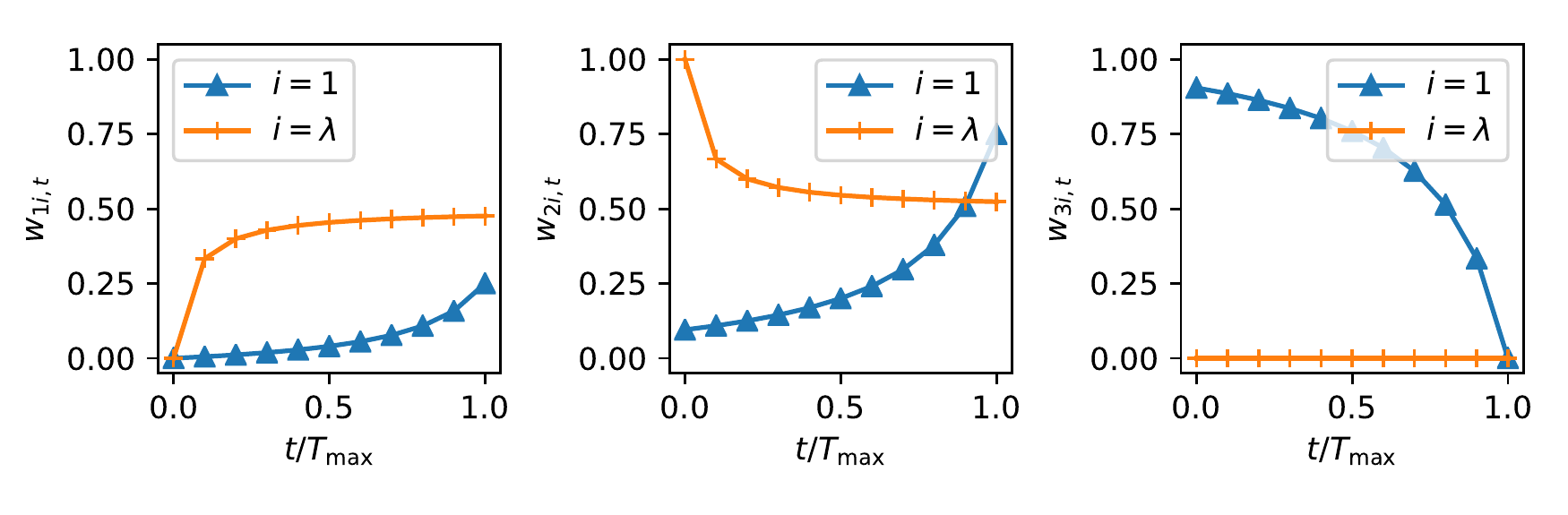}
 \caption{The change of weights for $1$st and $\lambda$th individuals on   CEC2017 benchmarks. $\gamma =0.1$.}
  \label{fig6} 
\end{figure}

Figures~\ref{fig5} and~\ref{fig6} depict the change of normalised weights  over  $t/T_{\max}$.  For $\lambda${th} individual, weights  $w_{1\lambda}>0, w_{2\lambda}>0$ but $w_{3\lambda}=0$. This individual minimises an equivalent function $f_{\lambda}$. For $1${st} individual, weight $w_{31}$ initially is set to a large value. Thus, at the beginning of search, this individual focuses on minimising a helper function $f_1$. Subsequently  $w_{31}$ decreases to $0$. It turns to minimise an equivalent function $f_1$ at the end of search.   

\section{Comparative Experiments and Results}
\label{secExperiments} 
\subsection{Experimental Setting}
HECO-DE was tested on two well-known benchmark sets. The first set is from IEEE  CEC2017 Competition and Special Session on Constrained Single Objective Real-Parameter Optimization~\cite{cec2017online} which consists of $28$ scalable functions with dimension $D=10, 30, 50, 100$ (total $4\times 28$ benchmarks).
The second set is from the IEEE CEC2006  Special Session on Constrained Real-parameter
Optimization~\cite{liang2006problem} which consists of 24 functions.  According to~\cite{liang2006problem}, there is no feasible solutions for function g20 and it is extremely difficult to find the optimum of function g22. Thus, these two functions are excluded in the comparison.

Tables~\ref{table:parameters1} and~\ref{table:parameters2} list the parameter setting used in HECO-DE. In Table~\ref{table:parameters1}, parameters inherited from LSHADE44 are set to    values similar to LSHADE44~\cite{polakova2017shade}.

\begin{table}[ht] \centering
        \caption{Parameters inherited from LSHADE44}
        \label{table:parameters1}
    \begin{tabular}{|c|c|}
    \hline 
    historical memory size & $H = 5$\\
        \hline
          number of strategies & $K = 4$ \\
       \hline
         constant in strategy adaption &$n_0 = 2$\\
       \hline
         threshold in strategy adaption & $\delta =1/20$\\
         \hline
         the maximum size of archive $A$ & $|A|_{\max} = 4N_t$\\ 
         \hline 
          tolerance for equivalent constraints & $\sigma = 0.0001$\\
 \hline
    \end{tabular}
\end{table} 

In Table~\ref{table:parameters2}, population size $N_0$, the number of subproblems $\lambda$ and constraint violation bias $\gamma$ are set to different values on CEC2006 and CEC2017 benchmarks. Since CEC2006 benchmarks include more constraints,  both the values of  $\lambda$ and   $\gamma$ are  set higher on CEC2006  benchmarks than that on CEC2017. The initial population size $N_0$ is set  to a constant on CEC2006 benchmarks, while it is set to $12D$ on CEC2017 benchmarks because the dimension $D$  ranges from 10 to 100.   As required by the competitions, twenty five independent runs were taken on each benchmark.

\begin{table}[ht] \centering
        \caption{Different parameter setting in CEC2006 and CEC2017}
        \label{table:parameters2}
    \begin{tabular}{|c|c|}
    \hline \multicolumn{2}{|c|}{CEC2006}\\\hline
       $FES_{\max}$ from CEC2006 benchmarks &  $FES_{\max}=500,000$\\ 
       \hline
                  required population sizes & $N_{0} = 450$,  $N_{T_{\max}} = \lambda$\\
                  \hline
          population size of $Q$ & $\lambda = 45$\\  \hline constraint violation bias  & $\gamma =0.7$\\ 
 \hline \multicolumn{2}{|c|}{CEC2017}\\ 
    \hline
     $FES_{\max}$ from CEC2017 benchmarks &  $FES_{\max}=20000 D$\\    \hline
           required population sizes & $N_{0} = 12D$,  $N_{T_{\max}} = \lambda$\\
       \hline
         population size of $Q$ & $\lambda = 20$\\   \hline constraint violation bias  & $\gamma =0.1$\\  
 \hline
    \end{tabular} 
\end{table}

\subsection{Experimental results on IEEE CEC2017 benchmarks}
HECO-DE was compared with seven single-objective EAs in CEC2017/18 constrained optimisation  competitions, which are CAL-SHADE~\cite{zamuda2017adaptive}, LSHADE44+IDE~\cite{tvrdik2017simple}, 
LSHADE44~\cite{polakova2017shade}, UDE~\cite{trivedi2017unified}, MA-ES~\cite{hellwig2018matrix}, IUDE~\cite{Trivedi2018improved}, LSHADE-IEpsilon~\cite{fan2018lshade44}, and one decomposition-based MOEA, DeCODE~\cite{wang2018decomposition}. 

HECO-DE was also compared with its two  variants. The first variant is to remove the equivalent function from HECO-DE. In the weighted sum (\ref{eq:fi}), $\tilde{e}(\vec{x})$ is replaced by $f(\vec{x})$. We call it HCO-DE.  The second variant is to choose the superiority of feasibility rule as the equivalent function.  In the weighted sum (\ref{eq:fi}), $\tilde{e}(\vec{x})$ is replaced by $e(\vec{x})$ given by (\ref{equFeasibleRule}). We call it HECO-DE(FR). The three algorithms adopt same parameter setting. 

According to the  CEC2017/18 competition rules~\cite{cec2017online},   EAs under comparison were ranked on the experimental results  against the use of 28 benchmarks under $D= 10, 30, 50, 100$, in terms of the mean values and median solution. All results were compared at the precision level of $1e-8$ in the same way as the official ranking source code~\cite{cec2017online}.  
The rank value of each algorithm on each dimension was 
calculated as below:
\begin{equation}
   \begin{aligned}
    \textrm{Rank value} &= \textstyle\sum_{i = 1}^{28}\mbox{rank}_i\textrm{(by mean value)}  \\
    &\textstyle+ \sum_{i = 1}^{28}\mbox{rank}_i\textrm{(by median solution)}.
   \end{aligned}
\end{equation}
The total rank value is the sum of rank values on four dimensions.

Table~\ref{table:rank1} summarises the ranks of EAs on  four dimensions and total ranks. HECO-DE is the top-ranked amongst all compared. This result clearly demonstrates that HECO-DE  consistently outperforms other EAs on all dimensions. Without the equivalent function, HCO-DE is worse than HECO-DE and HECO-DE(FR). HECO-DE(FR) which uses the  superiority of feasibility rule as the equivalent objective is slightly  worse than HECO-DE. Tables~\ref{table:rank2} and~\ref{table:rank3} provide a sensitivity analysis of parameters $\lambda$ and $\gamma$. HECO-DE with all five $\lambda$ and $\gamma$ values  had obtained lower total ranks than other EAs.

Due to the paper length restriction, more experimental results are provided in the supplement.

\begin{table}[ht]
\centering
\caption{Total ranks of HECO-DE and other EAs on IEEE CEC2017 benchmarks}
\label{table:rank1}
\begin{tabular}{|lccccc|}
\hline
Algorithm/Dimension    & $10D$        & $30D$        & $50D$        & $100D$       & Total        \\ \hline
CAL\_LSAHDE(2017)      & 421          & 420          & 469          & 478          & 1788         \\
LSHADE44+IDE(2017)     & 310          & 394          & 422          & 392          & 1518         \\
LSAHDE44(2017)         & 332          & 344          & 342          & 342          & 1360         \\
UDE(2017)              & 341          & 372          & 377          & 438          & 1528         \\
MA\_ES(2018)           & 271          & 261          & 273          & 282          & 1087         \\
IUDE(2018)             & 198          & 261          & 269          & 327          & 1055         \\
LSAHDE\_IEpsilon(2018) & 222          & 278          & 324          & 372          & 1196         \\
DeCODE(2018)           & 239          & 297          & 302          & 328          & 1166         \\
HCO-DE                 & 282          & 253          & 255          & 219          & 1009         \\
HECO-DE(FR)            & 158          & 194          & 186          & \textbf{202} & 740          \\
\textbf{HECO-DE}       & \textbf{154} & \textbf{139} & \textbf{156} & 205          & \textbf{654} \\ \hline
\end{tabular}
\end{table}

\begin{table}[ht]
\centering
\caption{Total ranks of HECO-DE with varying $\lambda$ and other EAs on IEEE CEC2017 benchmarks}
\label{table:rank2}
\begin{tabular}{|lccccc|}
\hline
Algorithm/Dimension            & $10D$        & $30D$        & $50D$        & $100D$       & Total        \\ \hline
CAL\_LSAHDE(2017)              & 507          & 508          & 569          & 582          & 2166         \\
LSHADE44+IDE(2017)             & 381          & 486          & 524          & 483          & 1874         \\
LSAHDE44(2017)                 & 409          & 431          & 431          & 422          & 1693         \\
UDE(2017)                      & 431          & 479          & 480          & 537          & 1927         \\
MA\_ES(2018)                   & 326          & 321          & 341          & 347          & 1335         \\
IUDE(2018)                     & 250          & 343          & 345          & 424          & 1362         \\
LSAHDE\_IEpsilon(2018)         & 277          & 354          & 420          & 472          & 1523         \\
DeCODE(2018)                   & 301          & 381          & 390          & 410          & 1482         \\
HECO-DE($\lambda=15$)          & \textbf{172} & 199          & 218          & 261          & 850          \\
\textbf{HECO-DE($\lambda=20$)} & 194          & \textbf{149} & \textbf{181} & 242          & \textbf{766} \\
HECO-DE($\lambda=25$)          & 177          & 174          & 197          & 241          & 789          \\
HECO-DE($\lambda=30$)          & 195          & 192          & 204          & \textbf{210} & 801          \\
HECO-DE($\lambda=35$)          & 189          & 208          & 200          & 222          & 819          \\ \hline
\end{tabular}
\end{table}

\begin{table}[ht]
\centering
\caption{Total ranks of HECO-DE with varying $\gamma$ values and other EAs on IEEE CEC2017 benchmarks}
\label{table:rank3}
\begin{tabular}{|lccccc|}
\hline
Algorithm/Dimension            & $10D$        & $30D$        & $50D$        & $100D$       & Total        \\ \hline
CAL\_LSAHDE(2017)              & 508          & 511          & 572          & 583          & 2174         \\
LSHADE44+IDE(2017)             & 373          & 485          & 518          & 482          & 1858         \\
LSAHDE44(2017)                 & 405          & 428          & 427          & 422          & 1682         \\
UDE(2017)                      & 423          & 471          & 465          & 532          & 1891         \\
MA\_ES(2018)                   & 329          & 320          & 334          & 349          & 1332         \\
IUDE(2018)                     & 249          & 317          & 315          & 419          & 1300         \\
LSAHDE\_IEpsilon(2018)         & 276          & 341          & 415          & 475          & 1507         \\
DeCODE(2018)                   & 296          & 362          & 370          & 398          & 1426         \\
HECO-DE($\gamma=0.0$)          & 254          & 207          & 243          & 287          & 991          \\
\textbf{HECO-DE($\gamma=0.1$)} & 186          & \textbf{177} & \textbf{186} & 234          & \textbf{783} \\
HECO-DE($\gamma=0.2$)          & \textbf{182} & 186          & 197          & 223          & 788          \\
HECO-DE($\gamma=0.3$)          & 190          & 220          & 229          & \textbf{210} & 849          \\
HECO-DE($\gamma=0.4$)          & 209          & 262          & 283          & 261          & 1015         \\ \hline
\end{tabular}
\end{table}

\subsection{Experimental results on IEEE CEC2006 benchmarks}
HECO-DE was compared with five EAs, which are CMODE~\cite{wang2012combining}, NSES~\cite{jiao2013novel}, FROFI~\cite{wang2015incorporating}, DW~\cite{peng2018novel} and DeCODE~\cite{wang2018decomposition}, on IEEE CEC2006 benchmarks. 

Table~\ref{table:cec2006results} summarises experiment results, where ``Mean'' and ``Std Dev'' denote the mean and standard deviation of objective function values, respectively. As suggested in~\cite{liang2006problem}, a successful run is a run during which an algorithm finds a feasible solution $\vec{x}$ satisfying $f(\vec{x}_{best}) - f(\vec{x}^*)\le 0.0001$, where $f(\vec{x}_{best})$ is the best solution found by the algorithm and  $f(\vec{x}^*)$ is the optimum. In Table~\ref{table:cec2006results}, ``*'' denotes that the algorithm satisfies this successful rule in 25 runs for a test problem.

As shown in Table~\ref{table:cec2006results}, the performance of HECO-DE is similar to NSES, FROFI, DeCODE, which can always find optimum of all test problems. HECO-DE performs better than CMODE and DW. CMODE cannot find the optimum of problem g21   and DW cannot find the optimum of g17 with 100\% success rate.

HECO-DE was also compared with HCO-DE and HECO-DE(FR) on four functions g02, g10, g21, and g23. Table~\ref{table:compare_eq_2006} shows that HECO-DE always find the optimum on all test functions. But without an equivalent objective, HCO-DE has a lower success rate or feasible rate. HECO-DE(FR)   faces performance degradation on g10, g21, and g23, probably because the superiority of feasibility rule has a higher selection pressure than the equivalent function~(\ref{equEquivalentFunction}).

\begin{table*}
\centering
\caption{Comparative experiment results on IEEE CEC2006 benchmarks. * denotes the number of satisfying successful rule}
\label{table:cec2006results}
\begin{adjustbox}{max width=\textwidth}
\begin{tabular}{|c|c|c|c|c|c|c|}
\hline
\begin{tabular}[c]{@{}c@{}} \\  \end{tabular} & \begin{tabular}[c]{@{}c@{}}CMODE\\ Mean$\pm$Std Dev\end{tabular} & \begin{tabular}[c]{@{}c@{}}NSES\\ Mean$\pm$Std Dev\end{tabular} & \begin{tabular}[c]{@{}c@{}}DW\\ Mean$\pm$Std Dev\end{tabular} & \begin{tabular}[c]{@{}c@{}}FROFI\\ Mean$\pm$Std Dev\end{tabular} & \begin{tabular}[c]{@{}c@{}}DeCODE\\ Mean$\pm$Std Dev\end{tabular} & \begin{tabular}[c]{@{}c@{}}HECO-DE\\ Mean$\pm$Std Dev\end{tabular} \\ \hline
g01                                                                                                                            & -1.5000E+01$\pm$0.00E+00*                                                                                                                    & -1.5000E+01$\pm$4.21E-30*                                                                                                                   & -1.5000E+01$\pm$5.02E-14*                                                                                                                 & -1.5000E+01$\pm$0.00E+00*                                                                                                                    & -1.5000E+01$\pm$0.00E+00*                                                                                                                     & -1.5000E+01$\pm$0.00E+00*                                                                                                                      \\ \hline
g02                                                                                                                            & -8.0362E-01$\pm$2.42E-08*                                                                                                                    & -8.0362E-01$\pm$2.41E-32*                                                                                                                   & -8.0362E-01$\pm$9.99E-08*                                                                                                                 & -8.0362E-01$\pm$1.78E-07*                                                                                                                    & -8.0362E-01$\pm$3.12E-09*                                                                                                                     & -8.0362E-01$\pm$1.21E-06*                                                                                                                      \\ \hline
g03                                                                                                                            & -1.0005E+00$\pm$5.29E-10*                                                                                                                    & -1.0005E+00$\pm$5.44E-19*                                                                                                                   & -1.0005E+00$\pm$4.27E-12*                                                                                                                 & -1.0005E+00$\pm$4.49E-16*                                                                                                                    & -1.0005E+00$\pm$4.00E-16*                                                                                                                     & -1.0005E+00$\pm$3.54E-09*                                                                                                                      \\ \hline
g04                                                                                                                            & -3.0666E+04$\pm$2.64E-26*                                                                                                                    & -3.0666E+04$\pm$2.22E-24*                                                                                                                   & -3.0666E+04$\pm$0.00E+00*                                                                                                                 & -3.0666E+04$\pm$3.71E-12*                                                                                                                    & -3.0666E+04$\pm$3.71E-12*                                                                                                                     & -3.0666E+04$\pm$0.00E+00*                                                                                                                      \\ \hline
g05                                                                                                                            & 5.1265E+03$\pm$1.24E-27*                                                                                                                     & 5.1265E+03$\pm$0.00E+00*                                                                                                                    & 5.1265E+03$\pm$4.22E-10*                                                                                                                  & 5.1265E+03$\pm$2.78E-12*                                                                                                                     & 5.1265E+03$\pm$2.78E-12*                                                                                                                      & 5.1265E+03$\pm$0.00E+00*                                                                                                                       \\ \hline
g06                                                                                                                            & -6.9618E+03$\pm$1.32E-26*                                                                                                                    & -6.9618E+03$\pm$0.00E+00*                                                                                                                   & -6.9618E+03$\pm$0.00E+00*                                                                                                                 & -6.9618E+03$\pm$0.00E+00*                                                                                                                    & -6.9618E+03$\pm$0.00E+00*                                                                                                                     & -6.9618E+03$\pm$0.00E+00*                                                                                                                      \\ \hline
g07                                                                                                                            & 2.4306E+01$\pm$7.65E-15*                                                                                                                     & 2.4306E+01$\pm$.37E-09*                                                                                                                     & 2.4306E+01$\pm$5.28E-10*                                                                                                                  & 2.4306E+01$\pm$6.32E-15*                                                                                                                     & 2.4306E+01$\pm$8.52E-12*                                                                                                                      & 2.4306E+01$\pm$1.77E-14*                                                                                                                       \\ \hline
g08                                                                                                                            & -9.5825E+02$\pm$6.36E-18*                                                                                                                    & -9.5825E+02$\pm$2.01E-34*                                                                                                                   & -9.5825E+02$\pm$2.78E-18*                                                                                                                 & -9.5825E+02$\pm$1.42E-17*                                                                                                                    & -9.5825E+02$\pm$1.42E-17*                                                                                                                     & -9.5825E+02$\pm$0.00E+00*                                                                                                                      \\ \hline
g09                                                                                                                            & 6.8063E+02$\pm$4.96E-14*                                                                                                                     & 6.8063E+02$\pm$1.10E-25*                                                                                                                    & 6.8063E+02$\pm$2.23E-11*                                                                                                                  & 6.8063E+02$\pm$2.23E-11*                                                                                                                     & 6.8063E+02$\pm$2.54E-13*                                                                                                                      & 6.8063E+02$\pm$5.57E-14*                                                                                                                       \\ \hline
g10                                                                                                                            & 7.0492E+03$\pm$2.52E-13*                                                                                                                     & 7.0492E+03$\pm$2.07E-24*                                                                                                                    & 7.0492E+03$\pm$4.43E-08*                                                                                                                  & 7.0492E+03$\pm$3.26E-12*                                                                                                                     & 7.0492E+03$\pm$6.34E-10*                                                                                                                      & 7.0492E+03$\pm$1.35E-06*                                                                                                                       \\ \hline
g11                                                                                                                            & 7.499E-01$\pm$0.00E+00*                                                                                                                      & 7.499E-01$\pm$0.00E+00*                                                                                                                     & 7.499E-01$\pm$1.06E-16*                                                                                                                   & 7.499E-01$\pm$1.13E-16*                                                                                                                      & 7.499E-01$\pm$1.13E-16*                                                                                                                       & 7.499E-01$\pm$0.00E+00*                                                                                                                        \\ \hline
g12                                                                                                                            & -1.00E+00$\pm$0.00E+00*                                                                                                                      & -1.00E+00$\pm$0.00E+00*                                                                                                                     & -1.00E+00$\pm$0.00E+00*                                                                                                                   & -1.00E+00$\pm$0.00E+00*                                                                                                                      & -1.00E+00$\pm$0.00E+00*                                                                                                                       & -1.00E+00$\pm$0.00E+00*                                                                                                                        \\ \hline
g13                                                                                                                            & 5.3942E-02$\pm$1.04E-17*                                                                                                                     & 5.3942E-02$\pm$1.98E-34*                                                                                                                    & 5.3942E-02$\pm$6.03E-14*                                                                                                                  & 5.3942E-02$\pm$2.41E-17*                                                                                                                     & 5.3942E-02$\pm$2.13E-17*                                                                                                                      & 5.3942E-02$\pm$1.30E-17*                                                                                                                       \\ \hline
g14                                                                                                                            & -4.7765E+01$\pm$3.62E-15*                                                                                                                    & -4.7765E+01$\pm$0.00E+00*                                                                                                                   & -4.7765E+01$\pm$3.47E-10*                                                                                                                 & -4.7765E+01$\pm$2.34E-14*                                                                                                                    & -4.7765E+01$\pm$2.93E-14*                                                                                                                     & -4.7765E+01$\pm$2.60E-15*                                                                                                                      \\ \hline
g15                                                                                                                            & 9.6172E+02$\pm$0.00E+00*                                                                                                                     & 9.6172E+02$\pm$0.00E+00*                                                                                                                    & 9.6172E+02$\pm$4.47E-13*                                                                                                                  & 9.6172E+02$\pm$5.80E-13*                                                                                                                     & 9.6172E+02$\pm$5.80E-13*                                                                                                                      & 9.6172E+02$\pm$0.00E+00*                                                                                                                       \\ \hline
g16                                                                                                                            & -1.9052E+00$\pm$2.64E-26*                                                                                                                    & -1.9052E+00$\pm$2.62E-30*                                                                                                                   & -1.9052E+00$\pm$0.00E+00*                                                                                                                 & -1.9052E+00$\pm$4.53E-16*                                                                                                                    & -1.9052E+00$\pm$4.53E-16*                                                                                                                     & -1.9052E+00$\pm$0.00E+00*                                                                                                                      \\ \hline
g17                                                                                                                            & 8.8535E+03$\pm$1.24E-27*                                                                                                                     & 8.8535E+03$\pm$2.51E-23*                                                                                                                    &\textbf{ 8.8802E+03$\pm$3.63E+01}                                                                                                                   & 8.8535E+03$\pm$0.00E+00*                                                                                                                     & 8.8535E+03$\pm$3.23E-08*                                                                                                                      & 8.8535E+03$\pm$2.98E-08*                                                                                                                       \\ \hline
g18                                                                                                                            & -8.6603E-01$\pm$6.51E-17*                                                                                                                    & -8.6603E-01$\pm$4.62E-33*                                                                                                                   & -8.6603E-01$\pm$3.30E-07*                                                                                                                 & -8.6603E-01$\pm$6.94E-16*                                                                                                                    & -8.6603E-01$\pm$2.47E-16*                                                                                                                     & -8.6603E-01$\pm$0.00E+00*                                                                                                                      \\ \hline
g19                                                                                                                            & 3.2656E+01$\pm$1.07E-10*                                                                                                                     & 3.2656E+01$\pm$1.52E-05*                                                                                                                    & 3.2656E+01$\pm$3.37E-07*                                                                                                                  & 3.2656E+01$\pm$2.18E-14*                                                                                                                     & 3.2656E+01$\pm$2.25E-14*                                                                                                                      & 3.2656E+01$\pm$4.17E-10*                                                                                                                       \\ \hline
g21                                                                                                                            & \textbf{2.6195E+01$\pm$5.34E+01}                                                                                                                      & 1.9372E+02$\pm$1.62E-22*                                                                                                                    & 1.9372E+02$\pm$3.66E-09*                                                                                                                  & 1.9372E+02$\pm$2.95E-11*                                                                                                                     & 1.9372E+02$\pm$4.82E-10*                                                                                                                      & 1.9372E+02$\pm$5.17E-11*                                                                                                                       \\ \hline
g23                                                                                                                            & -4.0006E+02$\pm$7.33E-11*                                                                                                                    & -4.0006E+02$\pm$9.08E-26*                                                                                                                   & -4.0006E+02$\pm$6.49E-06*                                                                                                                 & -4.0006E+02$\pm$1.71E-13*                                                                                                                    & -4.0006E+02$\pm$1.66E-05*                                                                                                                     & -4.0006E+02$\pm$4.37E-09*                                                                                                                      \\ \hline
g24                                                                                                                            & -5.5080E+00$\pm$.24E-28*                                                                                                                     & -5.5080E+00$\pm$0.00E+00*                                                                                                                   & -5.5080E+00$\pm$0.00E+00*                                                                                                                 & -5.5080E+00$\pm$9.06E-16*                                                                                                                    & -5.5080E+00$\pm$9.06E-16*                                                                                                                     & -5.5080E+00$\pm$0.00E+00*                                                                                                                      \\ \hline
*                                                          & \textbf{21}                                                          & \textbf{22}                                                         & \textbf{21}                                                       & \textbf{22}                                                          & \textbf{22}                                                           & \textbf{22}                                                            \\ \hline
\end{tabular}
\end{adjustbox}
\end{table*}

\begin{table}[ht]
\centering
\caption{Comparison of HECO-DE with HCO-DE and HECO-DE(FR) on functions g02, g10, g21, and g23}
\label{table:compare_eq_2006}
\begin{adjustbox}{max width= 0.47\textwidth}
\begin{tabular}{|c|c|c|c|}
\hline
\multirow{2}{*}{CEC2006} & \multicolumn{3}{c|}{Mean (Success Rate$\%$){[}Feasible Rate$\%${]}}                                             \\ \cline{2-4} 
                       & HCO-DE                         & HECO-DE(FR)                            & HECO-DE                                 \\ \hline
g02                    & -0.8032(96){[}100{]}   & \textbf{-0.8036(100){[}100{]}} & \textbf{-0.8036(100){[}100{]}}   \\ \hline
g10                    & 6815.3984(76){[}80{]}  & 7013.3762(80){[}84{]}          & \textbf{7049.2480(100){[}100{]}} \\ \hline
g21                    & 23.2469(12){[}12{]}    & 7.4898(4){[}4{]}               & \textbf{193.7245(100){[}100{]}}  \\ \hline
g23                    & -376.0544(96){[}100{]} & -376.0436(92){[}100{]}         & \textbf{-400.0551(100){[}100{]}} \\ \hline
\end{tabular}
\end{adjustbox}
\end{table}

\section{Conclusions}
\label{secConclusions}
This paper has proposed a helper and equivalent objective method for constrained optimisation.  
It is theoretically proven that  for a hard problem called  ``wide gap'',  using helper and equivalent objectives can  shorten the time of crossing the``wide gap''.   This general theoretical result shows the strengths of multi-objective EAs in solving COPs.

A case study has been conducted for validating our method. An algorithm, called HECO-DE, has been implemented which employs both  helper and equivalent objectives and reuses search operators from LSHADE44~\cite{polakova2017shade}. A new equivalent function and a new mechanism of dynamically weighting are designed in HECO-DE. Experimental results show that the overall performance of HECO-DE is ranked first when compared with other state-of-art EAs on CEC2017 benchmarks. HECO-DE also performs well on  CEC2006 benchmarks.  

For future work,  we will consider each constraint violation degree as an individual helper objective and then design a many helper and equivalent objectives EA for COPs.



\clearpage
\section*{Supplement: Experiments and Results}

This supplement provides  further details of the benchmark problems used for comparative experimental investigations and of experimental results and comparisons.
 
\subsection{Description of EAs under comparison on CEC2017 benchmarks}
The first seven EAs come from the CEC2017/18 constrained optimisation competitions~\cite{cec2017online}. The last one, DeCODE~\cite{wang2018decomposition}, was a  decomposition-based multi-objective EAs for constrained optimisation published in 2018.

\begin{enumerate}
    \item CAL-SHADE~\cite{zamuda2017adaptive}: Success-History based Adaptive Differential Evolution Algorithm including liner population size reduction, enhanced with adaptive constraint violation handling, i.e. adaptive $\epsilon$-constraint handling.
    \item LSHADE+IDE~\cite{tvrdik2017simple}: A simple framework for cooperation of two advanced adaptive DE variants. The search process is divided into two stages: (i) search feasible solutions via minimizing  the mean violation and stopped if a number of feasible solutions are found. (ii) minimize the function value until the stop condition is reached.
    \item LSHADE44~\cite{polakova2017shade}: Success-History based Adaptive Differential Evolution Algorithm including liner population size reduction, uses three different additional strategies compete, with the superiority of feasibility rule.
    \item UDE~\cite{trivedi2017unified}: Uses three trial vector generation strategies and two parameter settings. At each generation, UDE divides the current population into two sub-populations. In the first population, UDE employs all the three trial vector generation strategies on each target vector. For another one, UDE employs strategy adaption from learning experience from evolution in first population.
    \item MA-ES~\cite{hellwig2018matrix}: Combines the Matrix Adaptation Evolution Strategy for unconstrained optimization with well-known constraint handling techniques. It handles box-constraints by reflecting exceeding components into the predefined box. Additional in-/equality constraints are dealt with by application of two constraint handling techniques: $\epsilon$-level ordering and a repair step that is based on gradient approximation.
    \item IUDE~\cite{Trivedi2018improved}: An improved version of UDE. Different from UDE, local search and duplication operators have been removed, it employs a combination of $\epsilon$-constraint handling technique and the superiority of feasibility rule.
    \item LSHADE-IEpsilon~\cite{fan2018lshade44}: An improved $\epsilon$-constrained handling method (IEpsilon) for solving constrained single-objective optimization problems. The IEpsilon method adaptively adjusts the value of $\epsilon$ according to the proportion of feasible solutions in the current population. Furthermore, a new mutation operator DE/randr1*/1 is proposed.
    \item DeCODE~\cite{wang2018decomposition}: A recent decomposition-based EA made use of the weighted sum approach to decompose the transformed bi-objective problem into a number of scalar optimisation subproblems and then applied differential evolution to solve them. They designed a strategy of adjusting weights and a restart strategy to tackle COPs with complicated constraints.
\end{enumerate}

\subsection{The IEEE CEC2006 benchmark suit}
\begin{table}[ht]
	\centering
	\caption{Description of 24 Benchmark Functions from IEEE CEC2006
	where $D$ denotes dimension, $\rho$ the estimated ratio between the feasible area and the search space, $f(\vec{x}^*)$ the optimum objective function value}
	\label{table:benchmark}
	\resizebox{0.45\textwidth}{!}{
\begin{tabular}{|c|c|c|c|c|}
		\hline
		Function & D         & Type       & $\rho$    & $f(\mathbf{x}^*)$   \\ \hline
		g01      & 13        & Quadratic  & 0.0111\%  & -15.00000000      \\\hline
		g02      & 20        & Nonlinear  & 99.9971\% & -0.8036191041     \\\hline
		g03      & 10        & Polynomial & 0.0000\%  & -1.0005001000     \\\hline
		g04      & 5         & Quadratic  & 51.1230\% & -30665.5386717833 \\\hline
		g05      & 4         & Cubic      & 0.0000\%  & 5126.4967140071   \\\hline
		g06      & 2         & Cubic      & 0.0066\%  & -6961.8138755802  \\\hline
		g07      & 10        & Quadratic  & 0.0003\%  & 24.3062090682     \\\hline
		g08      & 2         & Nonlinear  & 0.8560\%  & -0.0958250414     \\\hline
		g09      & 7         & Polynomial & 0.5121\%  & 680.6300573744    \\\hline
		g10      & 8         & Linear     & 0.0010\%  & 7049.2480205287   \\\hline
		g11      & 2         & Quadratic  & 0.0000\%  & 0.7499000000      \\\hline
		g12      & 3         & Quadratic  & 4.7713\%  & -1.0000000000     \\\hline
		g13      & 5         & Nonlinear  & 0.0000\%  & 0.0539415140      \\\hline
		g14      & 10        & Nonlinear  & 0.0000\%  & -47.7648884595    \\\hline
		g15      & 3         & Quadratic  & 0.0000\%  & 961.7150222900    \\\hline
		g16      & 5         & Nonlinear  & 0.0204\%  & -1.9051552585     \\\hline
		g17      & 6         & Nonlinear  & 0.0000\%  & 8853.5338748065   \\\hline
		g18      & 9         & Quadratic  & 0.0000\%  & -0.8660254038     \\\hline
		g19      & 15        & Nonlinear  & 33.4761\% & 32.6555929502     \\\hline
		g20      & 24        & Linear     & 0.0000\%  & 0.2049794002      \\\hline
		g21      & 7         & Linear     & 0.0000\%  & 193.7245100697    \\\hline
		g22      & 22        & Linear     & 0.0000\%  & 236.4309755040    \\\hline
		g23      & 9         & Linear     & 0.0000\%  & -400.055100000    \\\hline
		g24      & 2         & Linear     & 79.6556\% & 5.5080132716      \\ \hline
	\end{tabular}
	}
\end{table}

\subsection{The IEEE CEC2017 Benchmark Suit}

IEEE CEC2017 benchmark suit is listed in Table~\ref{table:cec2017}  which consists of $4 \times 28$ problems with the dimension $D=10, 30, 50, 100$.

\begin{table}
\centering
\caption{Details of 28 Test Problems from IEEE CEC2018.   $I$ is the Number of Inequality Constraints, $E$ is the Number of Equality Constraints}

\label{table:cec2017}
\begin{scriptsize}
%
}                \\
                                                                               &                                                                                   &                                                                            &                                                                                     \\ \hline
\end{tabular}
\end{scriptsize}
\end{table}

\subsection{Convergence Speed of HECO-DE on IEEE CEC2006 Benchmark Suit}
\begin{figure}[ht]
\centering
  \includegraphics[width=0.48\textwidth]{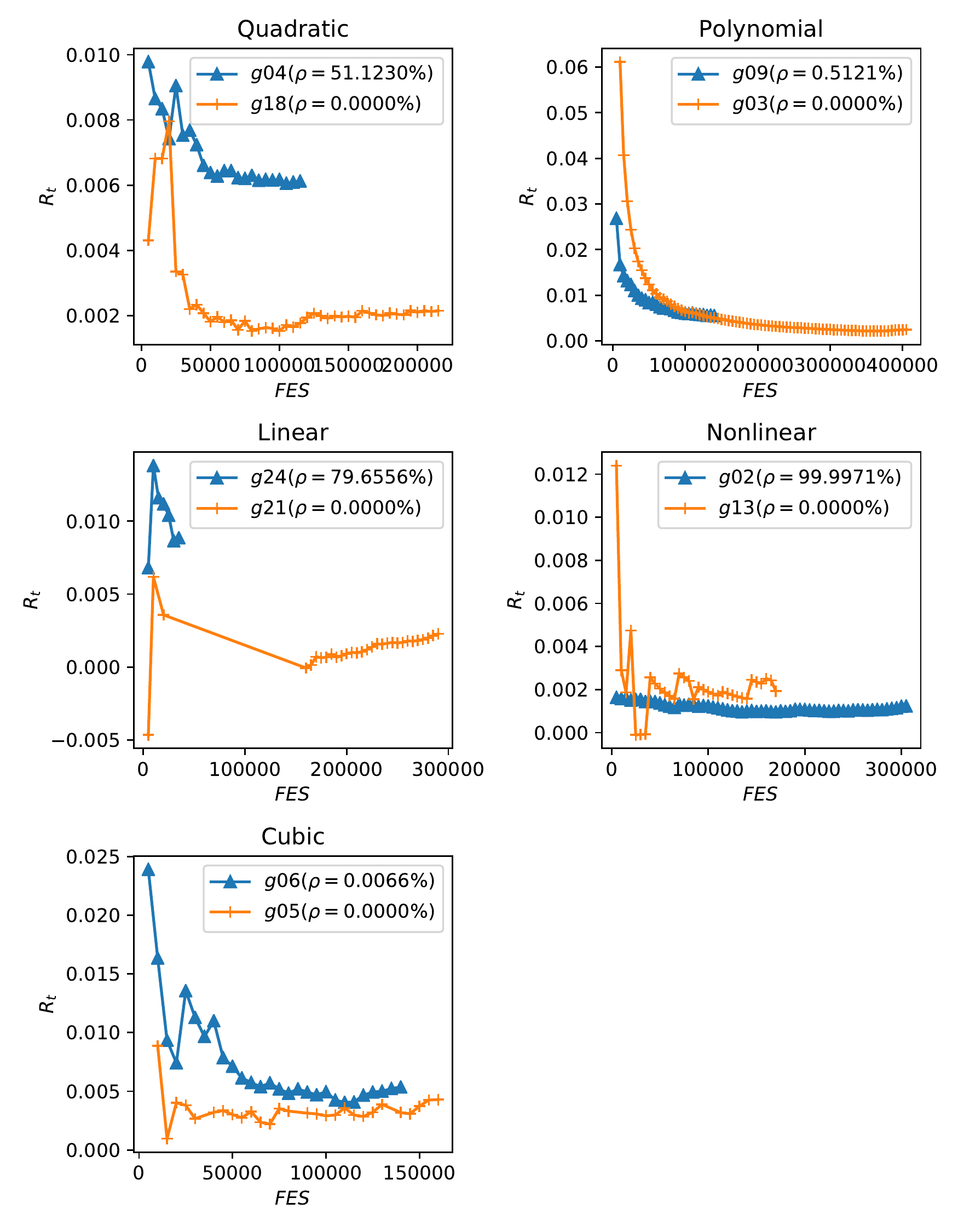}
 \caption{Average convergence rates on ten typical CEC2006 benchmark functions which are divided into five groups, such as quadratic, polynomial, linear, nonlinear and cubic, where $\rho$ denotes the estimated percentage of feasible region in the search space.}
  \label{fig7} 
\end{figure}

Fig.~\ref{fig7} plots the convergence speed at the median run of HECO-DE.
The convergence speed is measured by the average convergence rate $R_t$ defined as follows~\cite{he2016average}: 
\begin{equation}
R_t=1-\left|\frac{f_t-f^*}{f_0-f^*}\right|^{1/t}
\end{equation}
where $R_t$ denotes the normalised convergence speed, $t$ is the counter of the current generation, $f_t$ is the objective value at $t$ generation,   and $f^*$  the objective value of the known optimal solution.

Ten typical test function chosen from CEC2006 Benchmark are classified into five groups: quadratic, polynomial, linear, nonlinear and cubic. In each type of problems, we choose one function with relatively large feasible region and one function with very tiny feasible region.  

As shown in Fig.~\ref{fig7}, the convergence speed on all test functions is within the range around $[0.002, 0.01]$ after 50,000 generations.  The case of g12 is special. At the beginning, the convergence speed is negative. This implies an infeasible solution with $|f_t-f^*|>|f_0-f^*|$ is generated and accepted.

Fig.~\ref{fig7} shows that HECO-DE need more FES on test functions with tiny feasible region (g18, g03, g21 and g05) than test functions with large feasible region (g04, g09, g24 and g06) for satisfying the success criteria. However, this observation does not hold on nonlinear  functions (g13 and g0).

\subsection{Fine-tuning parameters on  CEC2006 benchmark}
CEC2006 benchmarks have more constraints than CEC2017 benchmarks. Thus the size of subpopulation $\lambda$ and constraint violation bias in CEC2006 are set to different values from CEC2017. 
Fine-tuning of parameters $\lambda$ and $\gamma$ was conducted on IEEE CEC2006 benchmark  functions.
For brevity, only performance on g02, g10, g17, g21, and g23 are shown in Tables~\ref{table:lambda2006} and~\ref{table:gamma2006} while other functions share the same performance with different value of parameter $\lambda$ and $\gamma$. As shown in Tables~\ref{table:lambda2006}, the value $\lambda = 45$ is the best because HECO-DE can always solve all tested benchmark functions 100\% successfully. As shown in Tables~\ref{table:gamma2006}, $\gamma = 0.7$ gives the best performance. The $\lambda$ and $\gamma$ values are larger than those used in CEC2017 ($\lambda=20$ and $\gamma =0.1$). This is due to CEC2006 benchmarks are strongly constrained.

\begin{table*}
\centering
\caption{Mean objective function value, success rate, feasible rate on IEEE CEC2006 benchmark functions g02, g10, g17, g21, and g23 with varied $\lambda$. }
\label{table:lambda2006}
\begin{tabular}{|c|c|c|c|c|c|}
\hline
\multirow{2}{*}{Prob.} & \multicolumn{5}{c|}{Mean (Success Rate$\%$){[}Feasible Rate$\%${]}}                                                                                                      \\ \cline{2-6} 
                       & 35                               & 40                               & 45                               & 50                               & 55                               \\ \hline
g02                    & -0.8032(92){[}100{]}             & \textbf{-0.8036(100){[}100{]}}   & \textbf{-0.8036(100){[}100{]}}   & \textbf{-0.8036(100){[}100{]}}   & -0.8032(96){[}100{]}             \\ \hline
g03                    & \textbf{-1.0005(100){[}100{]}}   & \textbf{-1.0005(100){[}100{]}}   & \textbf{-1.0005(100){[}100{]}}   & \textbf{-1.0005(100){[}100{]}}   & -1.00047(96){[}100{]}            \\ \hline
g10                    & \textbf{7049.2480(100){[}100{]}} & \textbf{7049.2480(100){[}100{]}} & \textbf{7049.2480(100){[}100{]}} & \textbf{7049.2480(100){[}100{]}} & 7049.2481(96){[}100{]}           \\ \hline
g13                    & \textbf{0.0539(100){[}100{]}}    & \textbf{0.0539(100){[}100{]}}    & \textbf{0.0539(100){[}100{]}}    & \textbf{0.0539(100){[}100{]}}    & 0.0539(96){[}100{]}              \\ \hline
g17                    & 8856.5008(96){[}100{]}           & \textbf{8853.5339(100){[}100{]}} & \textbf{8853.5339(100){[}100{]}} & \textbf{8853.5339(100){[}100{]}} & 8853.7232(96){[}100{]}           \\ \hline
g21                    & \textbf{193.7245(100){[}100{]}}  & \textbf{193.7245(100){[}100{]}}  & \textbf{193.7245(100){[}100{]}}  & \textbf{193.7245(100){[}100{]}}  & \textbf{193.7245(100){[}100{]}}  \\ \hline
g23                    & -388.0548(96){[}100{]}           & -376.0544(92){[}100{]}           & \textbf{-400.0551(100){[}100{]}} & -376.0544(92){[}100{]}           & \textbf{-400.0551(100){[}100{]}} \\ \hline
\end{tabular}
\end{table*}

\begin{table*}
\centering
\caption{Mean objective function value, success rate, feasible rate on IEEE CEC2006 benchmark functions g02, g10, g17, g21, and g23 with varied $\gamma$. }
\label{table:gamma2006}
\begin{tabular}{|c|c|c|c|c|c|}
\hline
\multirow{2}{*}{Prob.} & \multicolumn{5}{c|}{Mean (Success Rate$\%$){[}Feasible Rate$\%${]}}                                                                                                  \\ \cline{2-6} 
                       & 0.5                            & 0.6                            & 0.7                              & 0.8                              & 0.9                              \\ \hline
g02                    & \textbf{-0.8036(100){[}100{]}} & -0.8034(96){[}100{]}           & \textbf{-0.8036(100){[}100{]}}   & \textbf{-0.8036(100){[}100{]}}   & -0.8036(96){[}100{]}             \\ \hline
g03                    & \textbf{-1.0005(100){[}100{]}} & \textbf{-1.0005(100){[}100{]}} & \textbf{-1.0005(100){[}100{]}}   & \textbf{-1.0005(100){[}100{]}}   & \textbf{-1.0005(100){[}100{]}}   \\ \hline
g10                    & 10384.6108(8){[}92{]}          & 7049.2986(80){[}100{]}         & \textbf{7049.2480(100){[}100{]}} & \textbf{7049.2480(100){[}100{]}} & \textbf{7049.2480(100){[}100{]}} \\ \hline
g13                    & \textbf{0.0539(100){[}100{]}}  & \textbf{0.0539(100){[}100{]}}  & \textbf{0.0539(100){[}100{]}}    & \textbf{0.0539(100){[}100{]}}    & 0.0615(92){[}100{]}              \\ \hline
g17                    & 8854.9176(52){[}100{]}         & 8853.9032(80){[}100{]}         & \textbf{8853.5339(100){[}100{]}} & \textbf{8853.5339(100){[}100{]}} & 8856.5699(92){[}100{]}           \\ \hline
g21                    & 23.2469(12){[}12{]}            & 131.7327(68){[}68{]}           & \textbf{193.7245(100){[}100{]}}  & \textbf{193.7245(100){[}100{]}}  & \textbf{193.7245(100){[}100{]}}  \\ \hline
g23                    & -387.5537(80{[}100{]})         & -387.4869(88){[}100{]}         & \textbf{-400.0551(100){[}100{]}} & -376.0544(92){[}100{]}           & -376.0544(92){[}100{]}           \\ \hline
\end{tabular}
\end{table*}

\subsection{Detailed  experimental results  and ranking of HECO-DE on CEC2017 benchmarks}
In terms of IEEE CEC2017 benchmark functions, the best, median, worst, mean, standard deviation and feasibility rate of  the function values tested by HECO-DE on $10D$, $30D$, $50D$ and $100D$ are recorded in Table~\ref{table:fvalues10D}-\ref{table:fvalues100D}. 
\begin{itemize}
\item $c$ is the number of violated  constraints at the median solution where three figures indicate the number of violations (including inequality and equality) by more than 1.0, in the range
 $[0.01, 1.0]$ and in the range $[0.0001, 0.01]$ respectively. 
 
 \item $\overline{v}$ denotes the mean value of the constraint violations of all constraints at the median solution. 
 
 \item $SR$ is the 
 feasibility rate of the solutions obtained in 25 runs.
 \item $\overline{vio}$ denotes the mean
 constraint violation value of all the solutions in 25 runs.
 
\end{itemize}


As shown in Table~\ref{table:fvalues10D}-\ref{table:fvalues100D}, HECO-DE got high accuracy results with high feasibility rate on most test problems. However, no feasible solution was found in functions C17, C19, C26 and C28 on any dimensions. This is a common issue faced by all EAs when solving these problems. For functions C08, C11, C18, c22 and C27, a feasible solution sometimes was not found.

\subsection{Detailed ranking results of EAs on 2017 benchmarks}
For the 28 test problems in $10D$, $30D$, $50D$ and $100D$, the  ranks of each algorithm in terms of mean values and median solution are listed in Table~\ref{table:compare_mean_10D}-\ref{table:compare_median_100D} respectively. 

Regarding the test functions with $10D$, rank values based on mean values and median solution on the 28 test functions are reported in Table~\ref{table:compare_mean_10D} and~\ref{table:compare_median_10D}, respectively. In terms of mean of solutions, HECO-DE had the lowest rank values on $8$ of $28$ problems (functions C01-C03, C05-C09). However, HECO-DE got relatively poor performance on C11, C13, C16 and C25. HECO-DE got the second lowest total rank value $83$ which was slighter worse than the rank values obtained by HECO-DE(FR). In  terms of median solution, HECO-DE got the lowest rank value on $13$ of $28$ problems (functions C01-C09, C13, C16, C21, and C24). But its performance is not good on functions C11, C12 and C14. HECO-DE was ranked first with a total rank value $71$. The overall performance of HECO-DE is also the best among all nine EAs on $10D$ by summing up the two rank values in terms of mean values and median solution together. 

Regarding the test functions with $30D$, rank values based on mean values and median solution on the 28 test functions are listed in Table~\ref{table:compare_mean_30D} and~\ref{table:compare_median_30D}, respectively. 
HECO-DE had the lowest rank values on $11$ of $28$ problems (functions C01-C03, C06, C09, C10, C13, C15, 20, C21 and C24). However, HECO-DE got relatively poor performance on functions C05 and C11. In terms of median  solution, HECO-DE got the lowest rank value on $9$ of $28$ problems (functions C01-C03, C05, C06, C13, C15, C20 and C21). But its performance was not good on functions C11. Total rank values of HECO-DE were the lowest ones, $70$ in terms of mean of solutions and $69$ in terms of median solution, respectively.

Regarding the test functions with $50D$, rank values based on mean values and median solution on the 28 test functions are reported in Table~\ref{table:compare_mean_50D} and~\ref{table:compare_median_50D}, respectively. HECO-DE had the lowest rank values on $5$ of $28$ problems (functions C01-C05, C12, C15-C17, C21, C24 and C25). However, HECO-DE got relatively poor performance on functions C05 and C11. In terms of median solution, HECO-DE got the lowest rank value on $9$ of $28$ problems (functions  C01-C03, C05, C10, C12, C13, C20 and C23). But its performance was not good on functions C11. Total rank values of HECO-DE were the lowest ones, $88$ in terms of mean values and $68$ in terms of median solution, respectively.

Table~\ref{table:compare_mean_100D} and~\ref{table:compare_median_100D} record rank values based on mean values and median solution on the 28 test functions on $100D$. HECO-DE had the lowest rank values on $4$ of $28$ problems (functions C01, C02, C15 and C20). But HECO-DE got relatively poor performance on functions C05, C08, C11, C13 and C21. HECO-DE got the lowest total rank value $108$ here. In terms of median solution, HECO-DE got the lowest rank value on $5$ of $28$ problems (functions C01, C02, C15 and C20). But it had a  poor performance on functions C11-C13 and C21. HECO-DE got the second lowest total rank value $97$ which was only worse than the rank values obtained by HECO-DE(FR).

According to the competition rules, HECO-DE got the lowest or at least comparable total rank values on each dimension. This means that HECO-DE had an overall better performance than other eights algorithms on the IEEE CEC2017 benchmark suit. However, the ranking tables also show  that  no algorithm  could  perform better than other algorithms on all problems.

\begin{table*}
\centering
\caption{Function values of HECO-DE achieved for 10D ($FES_{\max} = 20000 \times D$) on IEEE CEC2017 benchmarks}
\label{table:fvalues10D}
\begin{scriptsize}
\begin{tabular}{|cccccccc|}
\hline
problem          & C01          & C02          & C03          & C04          & C05         & C06         & C07          \\ \hline
Best             & 0.00000e+00  & 0.00000e+00  & 0.00000e+00  & 0.00000e+00  & 0.00000e+00 & 0.00000e+00 & -9.91896e+02 \\
Median           & 0.00000e+00  & 0.00000e+00  & 0.00000e+00  & 0.00000e+00  & 0.00000e+00 & 0.00000e+00 & -9.68180e+02 \\
$c$              & 0 0 0        & 0 0 0        & 0 0 0        & 0 0 0        & 0 0 0       & 0 0 0       & 0 0 0        \\
$\overline{v}$   & 0.00000e+00  & 0.00000e+00  & 0.00000e+00  & 0.00000e+00  & 0.00000e+00 & 0.00000e+00 & 0.00000e+00  \\
mean             & 0.00000e+00  & 0.00000e+00  & 0.00000e+00  & 5.42911e+00  & 8.69719e-31 & 0.00000e+00 & -9.56102e+02 \\
Worst            & 0.00000e+00  & 0.00000e+00  & 0.00000e+00  & 1.35728e+01  & 2.17430e-29 & 0.00000e+00 & -8.80241e+02 \\
std              & 0.00000e+00  & 0.00000e+00  & 0.00000e+00  & 6.64928e+00  & 4.26074e-30 & 0.00000e+00 & 3.35462e+01  \\
\hline
\end{tabular}

\end{scriptsize}
\end{table*}

\begin{table*}
\centering
\caption{Ranks based on \textbf{mean solution} on the 28 functions of $10D$  on IEEE CEC2017 benchmarks}
\label{table:compare_mean_10D}
\scalebox{0.75}{
\begin{tabular}{|lccccccccccccccccccccccccccccc|}
\hline
\multicolumn{1}{|c}{Problem}                             & 1 & 2 & 3  & 4  & 5  & 6  & 7  & 8  & 9  & 10 & 11 & 12 & 13 & 14 & 15 & 16 & 17 & 18 & 19 & 20 & 21 & 22 & 23 & 24 & 25 & 26 & 27 & 28 & Total       \\ \hline
CAL\_LSAHDE(2017)                                        & 1 & 1 & 11 & 11 & 11 & 7  & 10 & 10 & 10 & 10 & 6  & 1  & 11 & 8  & 8  & 11 & 10 & 11 & 1  & 9  & 11 & 11 & 10 & 10 & 11 & 11 & 9  & 1  & 232         \\
LSHADE44+IDE(2017)                                       & 1 & 1 & 9  & 8  & 1  & 11 & 9  & 1  & 1  & 2  & 1  & 3  & 1  & 10 & 4  & 9  & 7  & 9  & 4  & 4  & 2  & 4  & 9  & 7  & 9  & 6  & 10 & 9  & 152         \\
LSAHDE44(2017)                                           & 1 & 1 & 10 & 6  & 1  & 10 & 8  & 1  & 9  & 2  & 2  & 10 & 1  & 9  & 9  & 10 & 8  & 8  & 2  & 1  & 7  & 8  & 8  & 6  & 10 & 7  & 8  & 10 & 173         \\
UDE(2017)                                                & 1 & 1 & 8  & 9  & 10 & 8  & 7  & 1  & 7  & 2  & 10 & 1  & 10 & 7  & 5  & 8  & 9  & 7  & 8  & 11 & 9  & 10 & 7  & 3  & 8  & 10 & 7  & 7  & 191         \\
MA\_ES(2018)                                             & 1 & 1 & 1  & 10 & 1  & 6  & 5  & 1  & 1  & 2  & 4  & 11 & 7  & 11 & 10 & 1  & 11 & 3  & 10 & 10 & 10 & 7  & 11 & 5  & 1  & 9  & 2  & 8  & 160         \\
IUDE(2018)                                               & 1 & 1 & 6  & 3  & 1  & 1  & 11 & 1  & 1  & 2  & 5  & 9  & 1  & 3  & 3  & 1  & 5  & 5  & 4  & 8  & 4  & 9  & 1  & 3  & 1  & 2  & 6  & 6  & 104         \\
LSAHDE\_IEpsilon(2018)                                   & 1 & 1 & 7  & 5  & 1  & 9  & 6  & 1  & 1  & 2  & 3  & 6  & 1  & 2  & 6  & 1  & 3  & 4  & 3  & 7  & 8  & 1  & 4  & 9  & 1  & 5  & 1  & 11 & 110         \\
\hline
\end{tabular}}
\end{table*} 

\begin{table*}
\centering
\caption{Ranks based on \textbf{median solution} on the 28 functions of $10D$ on IEEE CEC2017 benchmarks}
\label{table:compare_median_10D}
\scalebox{0.75}{
\begin{tabular}{|lccccccccccccccccccccccccccccc|}
\hline
\multicolumn{1}{|c}{Problem} & 1 & 2 & 3  & 4  & 5 & 6  & 7  & 8  & 9  & 10 & 11 & 12 & 13 & 14 & 15 & 16 & 17 & 18 & 19 & 20 & 21 & 22 & 23 & 24 & 25 & 26 & 27 & 28 & Total       \\ \hline
CAL\_LSAHDE(2017)            & 1 & 1 & 11 & 11 & 1 & 9  & 8  & 10 & 11 & 1  & 4  & 3  & 1  & 9  & 8  & 11 & 11 & 10 & 1  & 9  & 4  & 1  & 10 & 10 & 11 & 11 & 10 & 1  & 189         \\
LSHADE44+IDE(2017)           & 1 & 1 & 10 & 6  & 1 & 11 & 9  & 1  & 1  & 2  & 2  & 4  & 1  & 11 & 7  & 10 & 7  & 9  & 4  & 4  & 5  & 1  & 11 & 9  & 10 & 9  & 9  & 2  & 158         \\
LSAHDE44(2017)               & 1 & 1 & 9  & 8  & 1 & 10 & 10 & 1  & 1  & 2  & 4  & 11 & 1  & 10 & 9  & 9  & 8  & 7  & 2  & 1  & 1  & 1  & 9  & 8  & 9  & 7  & 8  & 10 & 159         \\
UDE(2017)                    & 1 & 1 & 7  & 9  & 1 & 8  & 7  & 1  & 9  & 2  & 10 & 1  & 1  & 1  & 5  & 8  & 10 & 8  & 10 & 11 & 1  & 1  & 1  & 5  & 8  & 10 & 7  & 6  & 150         \\
MA\_ES(2018)                 & 1 & 1 & 1  & 10 & 1 & 1  & 5  & 1  & 1  & 2  & 2  & 5  & 1  & 3  & 10 & 1  & 9  & 2  & 8  & 10 & 6  & 1  & 3  & 7  & 1  & 8  & 1  & 9  & 111         \\
IUDE(2018)                   & 1 & 1 & 1  & 1  & 1 & 1  & 11 & 1  & 1  & 2  & 1  & 10 & 1  & 3  & 5  & 1  & 1  & 6  & 4  & 8  & 9  & 1  & 4  & 5  & 1  & 1  & 6  & 6  & 94          \\
LSAHDE\_IEpsilon(2018)       & 1 & 1 & 8  & 5  & 1 & 1  & 6  & 1  & 1  & 2  & 6  & 8  & 1  & 3  & 3  & 1  & 6  & 3  & 3  & 7  & 11 & 1  & 8  & 4  & 1  & 6  & 2  & 11 & 112         \\
\hline
\end{tabular}}
\end{table*}

\begin{table*}
\centering
\caption{Function values of HECO-DE achieved for $30D$ ($FES_{\max} = 20000 \times D$)  on IEEE CEC2017 benchmarks}
\label{table:fvalues30D}
\begin{scriptsize}
\begin{tabular}{|cccccccc|}
\hline
problem          & C01          & C02          & C03          & C04          & C05         & C06         & C07          \\ \hline
Best             & 1.24862e-29  & 2.12623e-29  & 2.36658e-30  & 1.35728e+01  & 0.00000e+00 & 0.00000e+00 & -2.19162e+03 \\
Median           & 5.23113e-29  & 4.96613e-29  & 9.01274e-29  & 1.35728e+01  & 0.00000e+00 & 0.00000e+00 & -1.91485e+03 \\
$c$              & 0 0 0        & 0 0 0        & 0 0 0        & 0 0 0        & 0 0 0       & 0 0 0       & 0 0 0        \\
$\overline{v}$   & 0.00000e+00  & 0.00000e+00  & 0.00000e+00  & 0.00000e+00  & 0.00000e+00 & 0.00000e+00 & 0.00000e+00  \\
mean             & 5.63040e-29  & 5.57256e-29  & 1.07822e-28  & 1.35728e+01  & 1.59465e-01 & 4.92467e+01 & -1.90807e+03 \\
Worst            & 1.86036e-28  & 1.35548e-28  & 2.96612e-28  & 1.35728e+01  & 3.98662e+00 & 1.50462e+02 & -1.59126e+03 \\
std              & 3.48173e-29  & 2.56870e-29  & 7.57292e-29  & 5.65094e-15  & 7.81216e-01 & 6.11356e+01 & 1.59817e+02  \\
\hline
\end{tabular}
\end{scriptsize}
\end{table*}

\begin{table*}
\centering
\caption{Ranks based on \textbf{mean solution} on the 28 functions of $30D$  on IEEE CEC2017 benchmarks}
\label{table:compare_mean_30D}
\scalebox{0.75}{
\begin{tabular}{|lccccccccccccccccccccccccccccc|}
\hline
\multicolumn{1}{|c}{Problem} & 1 & 2 & 3  & 4  & 5  & 6  & 7  & 8  & 9  & 10 & 11 & 12 & 13 & 14 & 15 & 16 & 17 & 18 & 19 & 20 & 21 & 22 & 23 & 24 & 25 & 26 & 27 & 28 & Total       \\ \hline
CAL\_LSAHDE(2017)            & 1 & 1 & 11 & 10 & 11 & 5  & 11 & 1  & 11 & 10 & 4  & 11 & 11 & 11 & 9  & 11 & 11 & 10 & 1  & 6  & 11 & 11 & 11 & 10 & 11 & 11 & 9  & 1  & 232         \\
LSHADE44+IDE(2017)           & 1 & 1 & 10 & 6  & 1  & 11 & 7  & 10 & 1  & 9  & 3  & 8  & 8  & 10 & 7  & 8  & 8  & 9  & 4  & 8  & 9  & 8  & 10 & 9  & 8  & 8  & 10 & 9  & 201         \\
LSAHDE44(2017)               & 1 & 1 & 9  & 4  & 1  & 10 & 8  & 2  & 1  & 1  & 2  & 6  & 7  & 9  & 8  & 9  & 7  & 7  & 2  & 3  & 8  & 9  & 8  & 8  & 9  & 7  & 8  & 10 & 165         \\
UDE(2017)                    & 1 & 1 & 6  & 9  & 6  & 3  & 5  & 9  & 7  & 8  & 7  & 9  & 9  & 7  & 6  & 7  & 9  & 8  & 10 & 10 & 5  & 7  & 5  & 6  & 6  & 9  & 7  & 7  & 189         \\
MA\_ES(2018)                 & 1 & 1 & 1  & 8  & 1  & 2  & 4  & 2  & 1  & 1  & 1  & 10 & 1  & 8  & 10 & 1  & 10 & 2  & 11 & 11 & 10 & 1  & 7  & 7  & 1  & 10 & 1  & 5  & 129         \\
IUDE(2018)                   & 1 & 1 & 7  & 5  & 1  & 8  & 10 & 2  & 7  & 1  & 5  & 5  & 6  & 1  & 5  & 4  & 4  & 6  & 8  & 9  & 7  & 6  & 3  & 3  & 5  & 5  & 6  & 8  & 139         \\
LSAHDE\_IEpsilon(2018)       & 1 & 1 & 8  & 2  & 1  & 9  & 6  & 2  & 1  & 1  & 9  & 7  & 10 & 6  & 3  & 6  & 6  & 1  & 3  & 4  & 4  & 10 & 1  & 5  & 7  & 6  & 2  & 11 & 133         \\
\hline
\end{tabular}}
\end{table*}

\begin{table*}
\centering
\caption{Ranks based on \textbf{median solution} on the 28 functions of $30D$  on IEEE CEC2017 benchmarks}
\label{table:compare_median_30D}
\scalebox{0.75}{
\begin{tabular}{|lccccccccccccccccccccccccccccc|}
\hline
\multicolumn{1}{|c}{Problem} & 1 & 2 & 3  & 4  & 5 & 6  & 7  & 8  & 9  & 10 & 11 & 12 & 13 & 14 & 15 & 16 & 17 & 18 & 19 & 20 & 21 & 22 & 23 & 24 & 25 & 26 & 27 & 28 & Total       \\ \hline
CAL\_LSAHDE(2017)            & 1 & 1 & 11 & 11 & 1 & 7  & 11 & 1  & 1  & 1  & 4  & 9  & 11 & 8  & 9  & 11 & 11 & 9  & 1  & 6  & 4  & 11 & 7  & 10 & 11 & 11 & 8  & 1  & 188         \\
LSHADE44+IDE(2017)           & 1 & 1 & 10 & 3  & 1 & 11 & 9  & 10 & 2  & 10 & 3  & 2  & 1  & 11 & 7  & 9  & 9  & 10 & 4  & 8  & 10 & 8  & 10 & 9  & 9  & 6  & 10 & 9  & 193         \\
LSAHDE44(2017)               & 1 & 1 & 9  & 5  & 1 & 10 & 8  & 2  & 2  & 2  & 2  & 8  & 9  & 10 & 8  & 10 & 8  & 6  & 2  & 3  & 9  & 9  & 9  & 7  & 10 & 9  & 9  & 10 & 179         \\
UDE(2017)                    & 1 & 1 & 6  & 10 & 1 & 5  & 5  & 9  & 8  & 9  & 6  & 10 & 8  & 7  & 6  & 6  & 10 & 8  & 8  & 10 & 4  & 7  & 1  & 6  & 7  & 10 & 7  & 7  & 183         \\
MA\_ES(2018)                 & 1 & 1 & 1  & 9  & 1 & 4  & 4  & 2  & 2  & 2  & 1  & 11 & 1  & 9  & 10 & 1  & 7  & 1  & 11 & 11 & 11 & 1  & 8  & 8  & 1  & 7  & 1  & 5  & 132         \\
IUDE(2018)                   & 1 & 1 & 7  & 6  & 1 & 9  & 7  & 2  & 8  & 2  & 5  & 6  & 1  & 1  & 3  & 5  & 3  & 7  & 8  & 9  & 4  & 1  & 1  & 1  & 6  & 3  & 6  & 8  & 122         \\
LSAHDE\_IEpsilon(2018)       & 1 & 1 & 8  & 2  & 1 & 8  & 6  & 2  & 2  & 2  & 9  & 7  & 10 & 3  & 5  & 8  & 6  & 2  & 3  & 4  & 8  & 10 & 6  & 5  & 8  & 5  & 2  & 11 & 145         \\
\hline
\end{tabular}}
\end{table*}

\begin{table*}
\centering
\caption{Function Values of HECO-DE Achieved for $50D$ ($FES_{\max} = 20000 \times D$)  on IEEE CEC2017 benchmarks}
\label{table:fvalues50D}
\begin{scriptsize}
\begin{tabular}{|cccccccc|}
\hline
problem          & C01          & C02          & C03          & C04          & C05         & C06         & C07          \\ \hline
Best             & 3.26897e-28  & 2.74240e-28  & 6.07817e-28  & 1.35728e+01  & 1.99828e-28 & 1.37141e+02 & -3.27891e+03 \\
Median           & 7.56234e-28  & 6.20353e-28  & 2.16503e-27  & 1.35728e+01  & 1.07226e-27 & 3.18299e+02 & -2.70042e+03 \\
$c$              & 0 0 0        & 0 0 0        & 0 0 0        & 0 0 0        & 0 0 0       & 0 0 0       & 0 0 0        \\
$\overline{v}$   & 0.00000e+00  & 0.00000e+00  & 0.00000e+00  & 0.00000e+00  & 0.00000e+00 & 0.00000e+00 & 0.00000e+00  \\
mean             & 8.56437e-28  & 6.88503e-28  & 6.18894e+00  & 1.38003e+01  & 4.78395e-01 & 3.29350e+02 & -2.71356e+03 \\
Worst            & 3.02740e-27  & 1.76074e-27  & 7.74180e+01  & 1.69142e+01  & 3.98662e+00 & 4.76668e+02 & -1.74763e+03 \\
std              & 5.13174e-28  & 3.12836e-28  & 2.09877e+01  & 7.84266e-01  & 1.29550e+00 & 8.49272e+01 & 3.98328e+02  \\
\hline
\end{tabular}
\end{scriptsize}
\end{table*}

\begin{table*}
\centering
\caption{Ranks based on \textbf{mean solution} on the 28 functions of $50D$  on IEEE CEC2017 benchmarks}
\label{table:compare_mean_50D}
\scalebox{0.75}{
\begin{tabular}{|lccccccccccccccccccccccccccccc|}
\hline
\multicolumn{1}{|c}{Problem} & 1  & 2  & 3  & 4  & 5  & 6  & 7  & 8  & 9  & 10 & 11 & 12 & 13 & 14 & 15 & 16 & 17 & 18 & 19 & 20 & 21 & 22 & 23 & 24 & 25 & 26 & 27 & 28 & Total       \\ \hline
CAL\_LSAHDE(2017)            & 11 & 11 & 11 & 10 & 10 & 6  & 11 & 11 & 10 & 8  & 4  & 8  & 11 & 11 & 10 & 11 & 11 & 9  & 1  & 11 & 7  & 10 & 11 & 10 & 11 & 11 & 8  & 1  & 255         \\
LSHADE44+IDE(2017)           & 10 & 1  & 10 & 4  & 1  & 10 & 8  & 10 & 7  & 10 & 1  & 5  & 8  & 10 & 7  & 8  & 9  & 8  & 4  & 7  & 10 & 7  & 9  & 9  & 8  & 9  & 9  & 10 & 209         \\
LSAHDE44(2017)               & 1  & 1  & 9  & 1  & 1  & 9  & 7  & 1  & 3  & 1  & 2  & 10 & 7  & 9  & 8  & 9  & 8  & 7  & 3  & 3  & 11 & 8  & 7  & 8  & 7  & 8  & 7  & 9  & 165         \\
UDE(2017)                    & 1  & 1  & 6  & 9  & 11 & 5  & 5  & 9  & 1  & 9  & 5  & 7  & 10 & 7  & 6  & 6  & 10 & 10 & 10 & 8  & 4  & 9  & 3  & 6  & 5  & 10 & 10 & 6  & 189         \\
MA\_ES(2018)                 & 1  & 1  & 1  & 8  & 1  & 2  & 4  & 2  & 8  & 1  & 3  & 11 & 9  & 8  & 9  & 1  & 7  & 1  & 11 & 10 & 9  & 6  & 8  & 7  & 1  & 6  & 1  & 5  & 142         \\
IUDE(2018)                   & 1  & 1  & 7  & 7  & 1  & 11 & 10 & 7  & 1  & 1  & 8  & 4  & 6  & 4  & 3  & 5  & 4  & 6  & 9  & 9  & 5  & 5  & 2  & 3  & 3  & 3  & 6  & 8  & 140         \\
LSAHDE\_IEpsilon(2018)       & 1  & 1  & 8  & 2  & 7  & 8  & 6  & 8  & 6  & 7  & 9  & 9  & 5  & 6  & 5  & 7  & 6  & 2  & 7  & 5  & 6  & 11 & 1  & 5  & 6  & 7  & 2  & 11 & 164         \\
\hline
\end{tabular}}
\end{table*}

\begin{table*}
\centering
\caption{Ranks based on \textbf{median solution} on the 28 functions of $50D$  on IEEE CEC2017 benchmarks}
\label{table:compare_median_50D}
\scalebox{0.75}{
\begin{tabular}{|lccccccccccccccccccccccccccccc|}
\hline
\multicolumn{1}{|c}{Problem} & 1 & 2 & 3  & 4  & 5  & 6  & 7  & 8  & 9  & 10 & 11 & 12 & 13 & 14 & 15 & 16 & 17 & 18 & 19 & 20 & 21 & 22 & 23 & 24 & 25 & 26 & 27 & 28 & Total       \\ \hline
CAL\_LSAHDE(2017)            & 1 & 1 & 11 & 11 & 1  & 7  & 11 & 11 & 7  & 8  & 4  & 8  & 11 & 7  & 11 & 11 & 11 & 10 & 1  & 6  & 8  & 10 & 6  & 10 & 11 & 11 & 8  & 1  & 214         \\
LSHADE44+IDE(2017)           & 1 & 1 & 10 & 3  & 1  & 10 & 8  & 10 & 9  & 10 & 3  & 4  & 9  & 11 & 7  & 10 & 9  & 9  & 5  & 8  & 10 & 8  & 10 & 9  & 9  & 9  & 10 & 10 & 213         \\
LSAHDE44(2017)               & 1 & 1 & 9  & 5  & 1  & 9  & 7  & 1  & 3  & 1  & 1  & 11 & 8  & 10 & 9  & 9  & 8  & 7  & 4  & 3  & 11 & 9  & 9  & 8  & 8  & 8  & 7  & 9  & 177         \\
UDE(2017)                    & 1 & 1 & 6  & 10 & 11 & 6  & 5  & 9  & 1  & 9  & 5  & 6  & 10 & 8  & 6  & 7  & 10 & 8  & 11 & 9  & 1  & 7  & 3  & 7  & 6  & 10 & 9  & 6  & 188         \\
MA\_ES(2018)                 & 1 & 1 & 1  & 9  & 1  & 3  & 4  & 2  & 8  & 1  & 2  & 10 & 1  & 9  & 8  & 1  & 6  & 1  & 10 & 11 & 9  & 6  & 8  & 6  & 1  & 5  & 1  & 5  & 131         \\
IUDE(2018)                   & 1 & 1 & 7  & 8  & 1  & 11 & 10 & 5  & 1  & 1  & 7  & 5  & 1  & 5  & 1  & 6  & 4  & 6  & 2  & 10 & 5  & 5  & 3  & 1  & 4  & 4  & 6  & 8  & 129         \\
LSAHDE\_IEpsilon(2018)       & 1 & 1 & 8  & 2  & 1  & 8  & 6  & 8  & 3  & 7  & 9  & 9  & 7  & 1  & 5  & 8  & 7  & 2  & 8  & 5  & 6  & 11 & 5  & 5  & 7  & 7  & 2  & 11 & 160         \\
\hline
\end{tabular}}
\end{table*}

\begin{table*}
\centering
\caption{Function Values of HECO-DE Achieved for $100D$ ($FES_{\max} = 20000 \times D$)  on IEEE CEC2017 benchmarks}
\label{table:fvalues100D}
\begin{scriptsize}
\begin{tabular}{|cccccccc|}
\hline
problem          & C01         & C02         & C03          & C04          & C05         & C06         & C07          \\ \hline
Best             & 5.47366e-21 & 7.19883e-21 & 1.41057e+02  & 1.69142e+01  & 2.66253e-17 & 5.05460e+02 & -5.00027e+03 \\
Median           & 3.21897e-19 & 6.57495e-19 & 3.36356e+02  & 4.97476e+01  & 1.26216e-14 & 9.07974e+02 & -1.75646e+03 \\
$c$              & 0 0 0       & 0 0 0       & 0 0 0        & 0 0 0        & 0 0 0       & 0 0 0       & 0 0 0        \\
$\overline{v}$   & 0.00000e+00 & 0.00000e+00 & 0.00000e+00  & 0.00000e+00  & 0.00000e+00 & 0.00000e+00 & 0.00000e+00  \\
mean             & 7.51270e-19 & 4.62157e-18 & 3.24358e+02  & 6.02942e+01  & 1.59465e+00 & 8.88484e+02 & -2.04944e+03 \\
Worst            & 4.86463e-18 & 9.19137e-17 & 4.95920e+02  & 2.27844e+02  & 3.98662e+00 & 1.09726e+03 & -3.68566e+02 \\
std              & 1.05527e-18 & 1.78522e-17 & 9.03104e+01  & 4.07022e+01  & 1.95304e+00 & 1.11626e+02 & 1.17219e+03  \\
\hline
\end{tabular}

\end{scriptsize}
\end{table*}

\begin{table*}
\centering
\caption{Ranks based on \textbf{mean solution} on the 28 functions of $100D$  on IEEE CEC2017 benchmarks}
\label{table:compare_mean_100D}
\scalebox{0.75}{
\begin{tabular}{|lccccccccccccccccccccccccccccc|}
\hline
\multicolumn{1}{|c}{Problem} & 1  & 2  & 3  & 4  & 5  & 6  & 7  & 8  & 9  & 10 & 11 & 12 & 13 & 14 & 15 & 16 & 17 & 18 & 19 & 20 & 21 & 22 & 23 & 24 & 25 & 26 & 27 & 28 & Total        \\ \hline
CAL\_LSAHDE(2017)            & 10 & 10 & 11 & 11 & 6  & 5  & 11 & 11 & 7  & 9  & 5  & 8  & 10 & 11 & 10 & 11 & 11 & 10 & 1  & 11 & 11 & 8  & 11 & 11 & 11 & 11 & 8  & 1  & 251          \\
LSHADE44+IDE(2017)           & 11 & 11 & 9  & 2  & 2  & 10 & 7  & 6  & 2  & 7  & 1  & 2  & 5  & 10 & 5  & 8  & 9  & 8  & 3  & 7  & 8  & 7  & 7  & 8  & 8  & 9  & 9  & 8  & 189          \\
LSAHDE44(2017)               & 1  & 1  & 8  & 1  & 1  & 9  & 8  & 1  & 1  & 1  & 4  & 11 & 6  & 9  & 8  & 9  & 7  & 7  & 2  & 6  & 10 & 6  & 9  & 7  & 7  & 7  & 7  & 10 & 164          \\
UDE(2017)                    & 9  & 9  & 6  & 10 & 11 & 6  & 5  & 10 & 10 & 8  & 6  & 3  & 11 & 6  & 7  & 5  & 10 & 9  & 10 & 9  & 5  & 11 & 1  & 9  & 4  & 10 & 10 & 6  & 216          \\
MA\_ES(2018)                 & 1  & 1  & 1  & 8  & 10 & 2  & 3  & 2  & 9  & 4  & 2  & 10 & 2  & 8  & 9  & 1  & 6  & 1  & 11 & 10 & 9  & 4  & 8  & 5  & 1  & 6  & 2  & 5  & 141          \\
IUDE(2018)                   & 1  & 1  & 10 & 9  & 5  & 11 & 10 & 3  & 5  & 5  & 3  & 4  & 7  & 4  & 4  & 6  & 4  & 6  & 8  & 8  & 3  & 9  & 5  & 6  & 5  & 4  & 6  & 9  & 161          \\
LSAHDE\_IEpsilon(2018)       & 8  & 8  & 7  & 3  & 9  & 8  & 6  & 7  & 6  & 6  & 9  & 7  & 4  & 5  & 6  & 7  & 8  & 2  & 6  & 5  & 4  & 10 & 4  & 10 & 6  & 8  & 4  & 11 & 184          \\
\hline
\end{tabular}}
\end{table*}

\begin{table*}
\centering
\caption{Ranks based on \textbf{median solution} on the 28 functions of $100D$ on IEEE CEC2017 benchmarks}
\label{table:compare_median_100D}
\scalebox{0.75}{
\begin{tabular}{|lccccccccccccccccccccccccccccc|}
\hline
\multicolumn{1}{|c}{Problem} & 1  & 2  & 3  & 4  & 5  & 6  & 7  & 8  & 9  & 10 & 11 & 12 & 13 & 14 & 15 & 16 & 17 & 18 & 19 & 20 & 21 & 22 & 23 & 24 & 25 & 26 & 27 & 28 & Total       \\ \hline
CAL\_LSAHDE(2017)            & 10 & 10 & 11 & 11 & 8  & 6  & 11 & 11 & 7  & 9  & 5  & 7  & 10 & 3  & 10 & 11 & 11 & 10 & 1  & 3  & 4  & 8  & 6  & 11 & 11 & 11 & 10 & 1  & 227         \\
LSHADE44+IDE(2017)           & 11 & 11 & 10 & 1  & 4  & 9  & 8  & 7  & 2  & 8  & 4  & 3  & 5  & 11 & 6  & 9  & 8  & 7  & 3  & 8  & 10 & 7  & 8  & 9  & 10 & 8  & 7  & 9  & 203         \\
LSAHDE44(2017)               & 1  & 1  & 9  & 2  & 5  & 8  & 9  & 1  & 1  & 1  & 2  & 11 & 6  & 10 & 9  & 10 & 6  & 8  & 2  & 7  & 11 & 6  & 10 & 10 & 9  & 5  & 8  & 10 & 178         \\
UDE(2017)                    & 9  & 9  & 6  & 10 & 11 & 5  & 5  & 10 & 9  & 10 & 6  & 2  & 11 & 7  & 8  & 6  & 10 & 9  & 10 & 10 & 8  & 11 & 3  & 7  & 6  & 9  & 9  & 6  & 222         \\
MA\_ES(2018)                 & 1  & 1  & 1  & 8  & 10 & 3  & 3  & 2  & 11 & 5  & 1  & 9  & 3  & 9  & 4  & 1  & 5  & 1  & 11 & 11 & 9  & 5  & 9  & 5  & 1  & 6  & 1  & 5  & 141         \\
IUDE(2018)                   & 1  & 1  & 8  & 9  & 6  & 10 & 7  & 3  & 10 & 4  & 3  & 4  & 7  & 5  & 7  & 7  & 4  & 6  & 7  & 9  & 3  & 9  & 5  & 6  & 7  & 4  & 6  & 8  & 166         \\
LSAHDE\_IEpsilon(2018)       & 8  & 8  & 7  & 3  & 9  & 7  & 6  & 8  & 8  & 6  & 9  & 8  & 4  & 6  & 4  & 8  & 7  & 2  & 6  & 6  & 5  & 10 & 4  & 8  & 8  & 7  & 5  & 11 & 188         \\
\hline
\end{tabular}}
\end{table*}

\end{document}